[1994/12/01]% LaTeX date must December 1994 or later
\documentclass[11pt,letterpaper]{amsart}
\chardef\bslash=`\\ % p. 424, TeXbook
\usepackage{graphicx}
\usepackage[breaklinks=true]{hyperref}
\usepackage{mathtools}

\newtheorem{thm}{Theorem}[section]

\newtheorem{lem}[thm]{Lemma}
\newtheorem{prop}[thm]{Proposition}

\theoremstyle{definition}
\newtheorem{defn}[thm]{Definition}
\newtheorem{rem}[thm]{Remark}

\theoremstyle{remark}

\newcommand{\eval}[2][\right]{\relax
  \ifx#1\right\relax \left.\fi#2#1\rvert}

\begin{document}
\title{Strong convexity in flip-graphs}

\author[L. Pournin]{Lionel Pournin}
\address{Universit{\'e} Paris 13, Villetaneuse, France}
\email{lionel.pournin@univ-paris13.fr}

\author[Z. Wang]{Zili Wang}
\address{Sun Yat-Sen University, Shenzhen, Guangdong, China}
\email{wangzli6@mail.sysu.edu.cn}

%\date{Received on MONTH, YEAR}
%\issueinfo{VOL}{NUM}{MONTH}{YEAR}
%\doiinfo{10.1007/DOI-NUMBER}
\begin{abstract}
The set of triangulations of a surface $\Sigma$ with a prescribed set $X$ of vertices can be endowed with a graph structure $\mathcal{F}(\Sigma,X)$ called a flip-graph, whose edges connect two triangulations that differ by a single arc. It is known that when $X$ is the vertex set of a convex Euclidean polygon $\mathrm{P}$, the subgraph $\mathcal{F}_\varepsilon(\mathrm{P},X)$ induced in $\mathcal{F}(\mathrm{P},X)$ by the triangulations that contain a given arc $\varepsilon$ is strongly convex in the sense that all the geodesic paths in $\mathcal{F}(\mathrm{P},X)$ between two such triangulations remain in that subgraph. Here, we provide a related result that involves a triangle instead of an arc: we show that if the three edges of a triangle $\tau$ appear in (possibly distinct) triangulations along a geodesic path in $\mathcal{F}(\mathrm{P},X)$, then $\tau$ must belong to a triangulation in that path. More generally, we prove that certain $3$\nobreakdash-dimensional simplicial complexes related to the geodesics in $\mathcal{F}(\mathrm{P},X)$ are flag and provide two consequences. The first consequence is that $\mathcal{F}_\varepsilon(\mathrm{P},X)$ is not always strongly convex when $X$ is obtained from the vertex set of $P$ by adding just two points. The second, in the case when $\Sigma$ is a topological surface, is that the number of arc crossings between two triangulations does not allow to approximate their distance in $\mathcal{F}(\Sigma,X)$ by a factor of less than $3/2$.
\end{abstract}
\maketitle
%\tableofcontents

\section{Introduction}\label{PW.sec.0}

In order to study the properties of a surface, it is convenient to decompose it into elementary pieces such as triangles. This can be achieved by embedding a collection of pairwise non-crossing arcs in the surface, whose complement is a disjoint union of triangles. When the surface is equipped with a metric, we assume that these arcs are geodesics. Otherwise, they will be isotopy classes of curves. %In both of these cases the two extremities of an arc are allowed to coincide.
Such a collection of arcs is called a triangulation. Two triangulations that differ by a single arc can be thought of as related by a local operation called a \emph{flip} that removes an arc from one triangulation and replaces it with the only other arc such that the resulting set of arcs is a triangulation of the same surface. Observe that these triangulations share the same set of vertices as this operation does not change the set of the arcs' endpoints. One can therefore associate a graph to a surface $\Sigma$ equipped with a prescribed set $X$ of marked points. The vertices of this graph are the triangulations of $\Sigma$ that admit $X$ as their vertex set and its edges connect any two triangulations that are related by a flip. This graph, which we denote by $\mathcal{F}(\Sigma,X)$ from now on, is a \emph{flip-graph}.

Flip-graphs appear in a number of contexts in geometry and topology. Their simplest manifestation, when $\Sigma$ is a topological disk and $X$ a finite subset of its boundary, is a popular example. This particular flip-graph coincides with the flip-graph $\mathcal{F}(\mathrm{P},X)$ of a convex Euclidean polygon $\mathrm{P}$ where $X$ denotes the vertex set of $\mathrm{P}$. It is often studied in these terms and it is also the simplest manifestation of the flip-graph of a convex Euclidean polytope \cite{DeLoeraRambauSantos2010,GelfandZelevinskiiKapranov1990,Lee1989,Santos2000,Santos2005}. The geometry of $\mathcal{F}(\mathrm{P},X)$ is especially interesting because of its relation with the rotation distance between binary trees \cite{CulikWood1982,SleatorTarjanThurston1988} and has been investigated thoroughly~\cite{Fabila-MonroyFlores-PenalozaHuemerHurtadoUrrutiaWood2009,HurtadoNoy1999,ParlierPetri2018,Pournin2014,SleatorTarjanThurston1988}. The time-complexity of computing distances in $\mathcal{F}(\mathrm{P},X)$ has also attracted a lot of attention \cite{AichholzerMulzerPilz2015,ClearyMaio2018,Dorfer2026,LubiwPathak2015,Pilz2014}. Flip-graphs related to  $\mathcal{F}(\mathrm{P},X)$ can be obtained with objects other than triangulations \cite{CarrDevadoss2006,CeballosPilaud2016,FelsnerKleistMutzeSering2020,FominZelevinsky2003,HartungHoangMutzeWilliams2021} or by using a subset of the triangulations of $\mathrm{P}$ and a modified flip operation \cite{Pournin2017,Simion2003}.

Denote by $\mathcal{F}_\varepsilon(\Sigma,X)$ the subgraph induced in $\mathcal{F}(\Sigma,X)$ by the triangulations that contain a given arc $\varepsilon$. It is shown in \cite{SleatorTarjanThurston1988} that if $X$ is the vertex set of a convex Euclidean polygon $\mathrm{P}$, then $\mathcal{F}_\varepsilon(\mathrm{P},X)$ is \emph{strongly convex} in the sense that all the geodesic paths in $\mathcal{F}(\mathrm{P},X)$ between two vertices of $\mathcal{F}_\varepsilon(\mathrm{P},X)$ remain in that subgraph. This property is instrumental for the study of flip-graphs \cite{ClearyMaio2018,Pournin2014,SleatorTarjanThurston1988}. It is obtained by projecting the paths in $\mathcal{F}(\mathrm{P},X)$ to paths in $\mathcal{F}_\varepsilon(\mathrm{P},X)$. In a related combinatorial setting, the corresponding property was established in~\cite{CeballosPilaud2016,Williams2017} for the $1$-skeletons of all generalized associahedra that arise from the theory of cluster algebras~\cite{FominZelevinsky2003}. Our first main result is a similar property which involves a triangle instead of an arc.

\begin{thm}\label{PW.sec.0.thm.1}
Consider a convex Euclidean polygon $\mathrm{P}$ with vertex set $X$.~If the three edges of a triangle $\tau$ appear in triangulations along a geodesic path in $\mathcal{F}(\mathrm{P},X)$, then $\tau$ belongs to a triangulation in that path.
\end{thm}

This theorem is obtained by thinking of the paths in $\mathcal{F}(\mathrm{P},X)$ as certain $3$-dimensional triangulations which we call \emph{blow-up triangulations of $\mathrm{P}$}. %This space and these triangulations are geometric with respect to two dimensions and topological with respect to the third.
A blow-up triangulation is built, roughly, by gluing tetrahedra that correspond to the flips in the considered path. While blow-up triangulations are not triangulations in the usual sense and in particular they are not always simplicial complexes, the strong convexity result from~\cite{SleatorTarjanThurston1988} implies that a blow-up triangulation arising from a geodesic path in $\mathcal{F}(\mathrm{P},X)$ is a simplicial complex. We will establish Theorem \ref{PW.sec.0.thm.1} in the following form.

\begin{thm}\label{PW.sec.0.thm.0}
Consider a convex Euclidean polygon $\mathrm{P}$. If $X$ is the vertex set of $\mathrm{P}$, then all the blow-up triangulations of $\mathrm{P}$ that are built from geodesic paths in $\mathcal{F}(\mathrm{P},X)$ are flag simplicial complexes.
\end{thm}

The proof strategy is to project the paths in $\mathcal{F}(\mathrm{P},X)$ in the spirit of~\cite{SleatorTarjanThurston1988} except that the projection takes place with respect to a triangle within the blow-up triangulation that corresponds to these paths instead of with respect to an arc within each of the triangulations along the considered paths.

%Theorems \ref{PW.sec.0.thm.1} and \ref{PW.sec.0.thm.0} are powerful new tools to investigate the geometry of $\mathcal{F}(\mathrm{P},X)$.
We will also provide variants of Theorem \ref{PW.sec.0.thm.0} in the case when $X$ is a strict superset of the vertex set of $\mathrm{P}$. We will refer to the points in $X$ as \emph{flat vertices} when they belong to the interior of an edge of $\mathrm{P}$ and, by analogy with the non-Euclidean case, as \emph{punctures} when they belong to the interior of $\mathrm{P}$. It is well-known that, when $X$ contains six punctures or six flat vertices of $\mathrm{P}$, not only the subgraphs $\mathcal{F}_\varepsilon(\mathrm{P},X)$ are not always strongly convex but they are sometimes not even weakly convex in the sense that they contain none of the geodesic paths in $\mathcal{F}(\mathrm{P},X)$ between two of their vertices \cite{BoseHurtado2009,CeballosPilaud2016,DeLoeraRambauSantos2010,HurtadoNoyUrrutia1999}. Our second main result, proven using a variant of Theorem \ref{PW.sec.0.thm.0}, is that this remains the case down to two punctures or two flat vertices.

\begin{thm}\label{PW.sec.0.thm.2}
Consider two non-negative integers $i$ and $j$ that sum to at least $2$. There exists a convex polygon $\mathrm{P}$ and a set $X$ obtained by adding $i$ flat vertices and $j$ punctures to the vertex set of $\mathrm{P}$ such that, for some arc $\varepsilon$, the subgraph $\mathcal{F}_\varepsilon(\mathrm{P},X)$ is not weakly convex within $\mathcal{F}(\mathrm{P},X)$.
\end{thm}

As we shall see, the subgraph $\mathcal{F}_\varepsilon(\mathrm{P},X)$ is always strongly convex when $X$ is made of the vertices of $\mathrm{P}$ and just one flat vertex. In particular, Theorem~\ref{PW.sec.0.thm.2} is sharp for convex polygons with flat vertices.

Flip-graphs of orientable or non-orientable topological surfaces have been considered in \cite{BellDisarloTang2018,DisarloParlier2019,GultepeLeininger2017,PandaParlierPournin2025,ParlierPournin2017,ParlierPournin2018a,ParlierPournin2025} and the higher-dimensional analogues in~\cite{BjornerLutz2000,Pachner1991}. The strong convexity result from~\cite{SleatorTarjanThurston1988} in the case of a convex Euclidean polygon (or equivalently, of a topological disk) has been extended to arbitrary orientable topological surfaces in~\cite{DisarloParlier2019}. This result has been used to study the (coarse) geometry of the mapping class group of these surfaces~\cite{DisarloParlier2019} and to establish sharp bounds on the diameter of the quotient of $\mathcal{F}(\Sigma,X)$ by certain subgroups of the mapping class group~\cite{ParlierPournin2017,ParlierPournin2018a}.

Our third main result, an application of Theorem \ref{PW.sec.0.thm.0}, is on the problem of computing the distances within $\mathcal{F}(\Sigma,X)$ when $\Sigma$ is a topological surface. A popular procedure to estimate the distance $d(T^-,T^+)$ between two triangulations $T^-$ and $T^+$ in $\mathcal{F}(\Sigma,X)$ is based on the number of crossing arc pairs between these triangulations \cite{ClearyMaio2018,DisarloParlier2019}. Indeed, it is known that performing some flip in one of the triangulations makes the number of such crossings decrease~\cite{DisarloParlier2019} (this is also true for convex Euclidean punctured polygons \cite{HankeOttmannSchuierer1996}). As a consequence, one can build a path in $\mathcal{F}(\Sigma,X)$ between two triangulations $T^-$ and $T^+$ by greedily flipping arcs in one of them that make the number of arc crossings with the other decrease most. The length $\ell(T^-,T^+)$ of that path is then taken as an estimation of the desired distance. It is shown in \cite{ClearyMaio2018} that $\ell(T^-,T^+)$ is not always equal to $d(T^-,T^+)$. Here we show that this estimation is sometimes off by a factor of about $3/2$.

\begin{thm}\label{PW.sec.0.thm.3}
Consider any positive number $\epsilon$ and any topological surface $\Sigma$. For any sufficiently large $n$, one can equip $\Sigma$ with a set $X$ of $n$ marked points in such a way that there exist two triangulations $T^-$ and $T^+$ whose distance in $\mathcal{F}(\Sigma,X)$ goes to infinity with $n$ and
$$
\frac{\ell(T^-,T^+)}{d(T^-,T^+)}\geq\frac{3}{2}-\epsilon\mbox{.}
$$
\end{thm}

This follows, via Theorem \ref{PW.sec.0.thm.0} and the decomposition strategy described in Section \ref{PW.sec.1}, from determining the exact distance in $\mathcal{F}(\Sigma,X)$ of a family of triangulation pairs for which that computation was so far out of reach.

The outline of the article is the following. In Section \ref{PW.sec.1}, we explain how paths in the flip-graph of a convex polygon, possibly with flat vertices can be modeled as blow-up triangulations. We derive some of the geometric and combinatorial properties of these triangulations in Section~\ref{PW.sec.1.1}. In particular, we translate in terms of blow-up triangulations the tools from \cite{Pournin2014} that will allow us to bound distances in flip-graphs. In Section \ref{PW.sec.1.5}, we recall how the projection from \cite{SleatorTarjanThurston1988} with respect to an arc in the flip-graph of a convex polygon works, and show that it still works for a convex polygon with a single flat vertex. In Section \ref{PW.sec.2}, we introduce our projection with respect to triangles and prove a decomposition lemma for blow-up triangulations. In Section \ref{PW.sec.2.5}, we give several variants of that decomposition lemma, one of which is Theorem \ref{PW.sec.0.thm.1} and another Theorem \ref{PW.sec.0.thm.0}. We establish Theorem~\ref{PW.sec.0.thm.2} in Section \ref{PW.sec.3} using a variant of a the decomposition lemma from Section \ref{PW.sec.2.5} that deals with convex polygons with flat vertices. Finally, in Section \ref{PW.sec.4}, we prove Theorem \ref{PW.sec.0.thm.3} as a consequence of Theorem \ref{PW.sec.0.thm.0}.

\section{A topological model for sequences of flips}\label{PW.sec.1}

As remarked in \cite{SleatorTarjanThurston1988}, a sequence of flips between two triangulations of a convex polygon can be thought of as a $3$-dimensional triangulation. Before we give a precise description of this correspondence in the more general context of convex polygons with flat vertices, let us illustrate it informally and give a quick word on how we will use these triangulations.

\begin{figure}[b]
\begin{centering}
\includegraphics[scale=1]{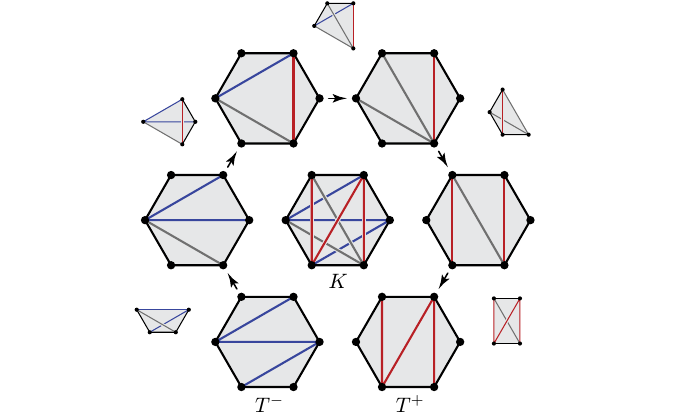}
\caption{A path between two of triangulations $T^-$ and $T^+$ of a convex hexagon and the corresponding blow-up triangulation $K$ of that hexagon (center). The tetrahedron in $K$ that corresponds to each flip is depicted next to it.}\label{PW.sec.1.fig.0}
\end{centering}
\end{figure}

Consider the path in $\mathcal{F}(\mathrm{P},X)$ shown in Fig. \ref{PW.sec.1.fig.0} between two triangulations $T^-$ and $T^+$, where $\mathrm{P}$ is a convex hexagon and $X$ its vertex set. Each flip in that path removes a diagonal of a convex quadrilateral from a triangulation and replaces it with the other diagonal of that quadrilateral. The tetrahedron depicted next to each flip is obtained from that quadrilateral by raising above the plane one of the vertices of the arc introduced by the flip while keeping the other three vertices in the plane. Now consider the $3$-dimensional triangulation $K$ obtained by gluing these tetrahedra on top of one another in the order of the flips they correspond to. This triangulation, sketched at the center of Fig. \ref{PW.sec.1.fig.0}, is an example of what we call a blow-up triangulation of $\mathrm{P}$. Note that this construction is topological with respect to the third dimension (so the tetrahedra can be bent up and down) and geometric with respect to the first two dimensions (in particular, the orthogonal projection of each tetrahedron back on the plane is the quadrilateral whose diagonals are exchanged by the corresponding flip).

\begin{figure}[b]
\begin{centering}
\includegraphics[scale=1]{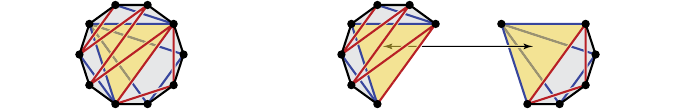}
\caption{When a blow-up triangulation $K$ of a convex polygon (left) contains a triangle (colored yellow) whose three edges lie in the boundary of $K$, one can decompose $K$ into two blow-up triangulations of smaller polygons (right).}\label{PW.sec.1.fig.1}
\end{centering}
\end{figure}

Observe that if a blow-up triangulation $K$ of a convex polygon $\mathrm{P}$ contains a triangle whose three edges are in the boundary of $K$, then we can split $K$ along that triangle, as sketched in Fig. \ref{PW.sec.1.fig.1}. This results into two blow-up triangulations of smaller polygons. More generally, $K$ can be cut along the union of some of its triangles provided that the boundary of that union is contained in the boundary of $K$ (see for instance Figures \ref{PW.sec.3.fig.5} and~\ref{PW.sec.4.fig.4}). The path in $\mathcal{F}(\mathrm{P},X)$ corresponding to $K$ can be then be analyzed by looking at smaller blow-up triangulations. Theorem~\ref{PW.sec.0.thm.0} and its variants established in Section \ref{PW.sec.2.5} provide the existence of triangles along which $K$ can be cut. We will use these cuts to prove Theorems \ref{PW.sec.0.thm.2} and \ref{PW.sec.0.thm.3}.

We now give a detailed description of blow-up triangulations. In the remainder of this section and in the next four sections, $\mathrm{P}$ is a convex Euclidean polygon contained in $\mathbb{R}^2$ and $X$ is a set that contains all the vertices of $\mathrm{P}$ together with finitely-many flat vertices. Consider a triangulation $T$ of $P$ with vertex set $X$ or equivalently, a vertex of $\mathcal{F}(\mathrm{P},X)$. The arcs of $T$ are Euclidean line segment between two points of $X$ that do not contain any other point of $X$ than their extremities. More precisely, $T$ is a set of pairwise non-crossing such arcs that is maximal with respect to inclusion. A triangle whose three edges belong to $T$ will be called a triangle \emph{of} $T$.

We identify $\mathbb{R}^2$ with the plane spanned by the first two coordinates of $\mathbb{R}^3$ and denote by $\pi$ the orthogonal projection from $\mathbb{R}^3$ to $\mathbb{R}^2$.

\begin{defn}\label{PW.sec.1.defn.1}
Given two subsets $\sigma$ and $\tau$ of $\mathbb{R}^3$, we say that $\sigma$ is \emph{below} $\tau$ (and that $\tau$ is \emph{above} $\sigma$) when $\pi(\sigma)$ is non-disjoint from $\pi(\tau)$ and, whenever a point $x$ of $\sigma$ has the same image by $\pi$ than a point $y$ of $\tau$, the third coordinate of $x$ is less than the third coordinate of $y$.
\end{defn}

%While definition \ref{PW.sec.1.defn.1} is stated for any subsets $\sigma$ and $\tau$ of $\mathbb{R}^3$, we will exclusively use it in the case when $\sigma$ and $\tau$ are topological balls.
Given two points $x$ and $y$ in $\mathbb{R}^3$, we shall say for short that $x$ is below $y$ (and $y$ above $x$) when $\{x\}$ is below $\{x\}$ in the sense of Definition \ref{PW.sec.1.defn.1}.

Throughout the article, we identify a path in $\mathcal{F}(\mathrm{P},X)$ between two triangulations $T^-$ and $T^+$ with the sequence $(T_i)_{0\leq{i}\leq{k}}$ of the triangulations in that path, ordered in such a way that $T^-$ is equal to $T_0$, $T^+$ to $T_k$, and $T_i$ is the triangulation that results from the $i$-th flip along the path (in particular, $T_{i-1}$ and $T_i$ are related by a flip). Let us explain how a blow-up triangulation of $\mathrm{P}$ can be obtained from a path $(T_i)_{0\leq{i}\leq{k}}$ in $\mathcal{F}(\mathrm{P},X)$. First consider a topological disk $\Pi_0$ embedded in $\mathbb{R}^3$ in such a way that $\pi$ induces a homeomorphism $\Pi_0\rightarrow\mathrm{P}$ and the boundaries of $\Pi_0$ and $\mathrm{P}$ coincide. If $k$ is equal to $0$, then the set made of the portions of $\Pi_0$ whose images by $\pi$ are a vertex, an arc, or a triangle of $T_0$ is the desired blow-up triangulation of $\mathrm{P}$. If however $k$ is positive, consider the portion $\Theta_1^-$ of $\Pi_0$ such that $\pi(\Theta_1^-)$ is the quadrilateral whose diagonals are exchanged by the flip between $T_0$ and $T_1$. We embed a topological disk $\Theta_1^+$ in $\mathbb{R}^3$ in such a way that $\pi$ induces a homeomorphism $\Theta_1^+\rightarrow\Theta_1^-$ while the boundary of $\Theta_1^+$ coincides with the boundary of $\Theta_1^-$ and every point in the interior of $\Theta_1^+$ is above the point from $\Pi_0$ with the same image by $\pi$. Denote
$$
\Pi_1=[\Pi_0\mathord{\setminus}\Theta_1^-]\cup\Theta_1^+
$$
and observe that $\pi$ induces a homeomorphism $\Pi_1\rightarrow\mathrm{P}$. Further observe that  $\Theta_1^+\cup\Theta_1^-$ is a $2$-dimensional topological sphere and denote by $\psi(1)$ the closed $3$-dimensional topological ball bounded by that sphere.

If $k$ is equal to $1$, the desired blow-up triangulation of $\mathrm{P}$ is the set made of $\psi(1)$, of the portions of $\Pi_0$ whose images by $\pi$ are a vertex, an arc, or a triangle of $T_0$, and of the portions of $\Pi_1$ whose images by $\pi$ are a vertex, an arc, or a triangle of $T_1$. Note that $\Theta_1^+$ is then decomposed into two triangles of $K$ sharing an arc whose image by $\pi$ is the arc of $T_1$ introduced by the first flip along the considered path. If $k$ is greater than $1$, one can repeat this procedure in order to build three families $(\Pi_i)_{0\leq{i}\leq{k}}$, $(\Theta_i^-)_{1\leq{i}\leq{k}}$, and $(\Theta_i^+)_{1\leq{i}\leq{k}}$ of topological disks that satisfy the following statement.

\begin{prop}
For any integer $i$ satisfying  $0\leq{i}\leq{k}$, $\pi$ induces a homeomorphism $\Pi_i\rightarrow\mathrm{P}$ and if $i$ is positive, then
\begin{itemize}
\item[(i)] $\Pi_i=[\Pi_{i-1}\mathord{\setminus}\Theta_i^-]\cup\Theta_i^+$,
\item[(ii)] all the points in the interior of $\Theta_i^+$ are above $\Pi_{i-1}$, and
\item[(iii)] the image by $\pi$ of both $\Theta_i^-$ and $\Theta_i^+$ is the quadrilateral whose diagonals are exchanged by the flip that transforms $T_{i-1}$ into $T_i$.
\end{itemize}
\end{prop}

As above, $\Theta_i^-\cup\Theta_i^+$ is a $2$-dimensional topological sphere that bounds a $3$-dimensional closed topological ball which we denote by $\psi(i)$. Just as for $\Theta_i^-$ and $\Theta_i^+$, the image by $\pi$ of $\psi(i)$ is the quadrilateral whose diagonals are exchanged by the flip that transforms $T_{i-1}$ into $T_i$. Now consider the set $K$ made of the portions of $\Pi_i$ whose image by $\pi$ is a vertex, an arc, or a triangle of $T_i$, when $i$ ranges from $0$ to $k$ together with the balls $\psi(i)$ when $i$ ranges from $1$ to $k$. Observe that all the elements of $K$ are then closed topological balls whose dimension ranges from $0$ to $3$.

\begin{defn}
A blow-up triangulation of $\mathrm{P}$ is a set of topological balls built from a path in $\mathcal{F}(\mathrm{P},X)$ just as $K$ is built from $(T_i)_{0\leq{i}\leq{k}}$.
\end{defn}

We will refer to the points contained in $K$ as its vertices. Note that by construction, the vertices of $K$ are exactly the elements of $X$. We will also call the $1$-, $2$-, and $3$-dimensional topological balls contained in $K$ as its arcs, triangles, and tetrahedra respectively. Note that all the arcs $\varepsilon$ contained in the boundary of $\mathrm{P}$ whose intersection with $X$ coincides with their vertex set are contained in any blow-up triangulation of $\mathrm{P}$ with vertex set $X$. In the sequel, we call the vertices of such an arc \emph{consecutive} in $X$. The triangles of $K$ are really topological disks whose boundary is the union of three arcs of $K$. Likewise, the tetrahedra of $K$ are really just $3$-dimensional topological balls but they are tetrahedra with respect to $K$ in the following sense.

\begin{rem}
If $\sigma$ is a tetrahedron of $K$, then there exists a homeomorphism $\phi$ from $\sigma$ to a Euclidean tetrahedron that sends the set of the elements of $K$ contained in $\sigma$ to the face complex of $\phi(\sigma)$.
\end{rem}

We will alternatively call the elements of $K$ its \emph{faces}. For any face $\sigma$ of $K$, the faces of $K$ contained in $\sigma$ will also be referred to as the faces of $\sigma$. The $0$-dimensional and $1$-dimensional faces of $\sigma$ will be called its vertices and its edges, respectively. A face of $\sigma$ whose dimension is less by one than the dimension of $\sigma$ will be referred to as a facet of $\sigma$.

Observe that a face $\sigma$ of $K$ cannot be both below and above another face $\tau$ of $K$. Indeed, by continuity and by the convexity of their image by $\pi$, the interiors of $\sigma$ and $\tau$ would otherwise not be disjoint.

\begin{rem}\label{PW.sec.1.rem.2}
Blow-up triangulations are not triangulations in the usual sense as they are not always simplicial complexes. For instance, two arcs of a blow-up triangulation $K$ can have the same pair of endpoints.
%Likewise, two tetrahedra can meet along the union of two facets
This happens for instance when $K$ is built from a path in $\mathcal{F}(\mathrm{P},X)$ along which a flip is immediately followed by the inverse flip. %Still, the intersection of two distinct faces of a blow-up triangulation $K$ is contained in the boundary of at least one of them and in turn, distinct faces of $K$ have disjoint interiors. In particular, a triangle of $K$ cannot be a face of more than two tetrahedra of $K$.
\end{rem}

\begin{figure}[b]
\begin{centering}
\includegraphics[scale=1]{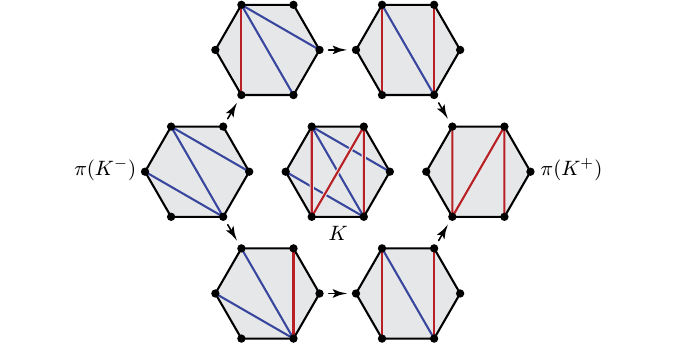}
\caption{A blow-up triangulation $K$ of a hexagon shown from above (center) and two different paths from $\pi(K^-)$ to $\pi(K^+)$ in associated to it (top and bottom).}\label{PW.sec.1.fig.2}
\end{centering}
\end{figure}

While different blow-up triangulations of $\mathrm{P}$ can be built from a path $(T_i)_{0\leq{i}\leq{k}}$ in $\mathcal{F}(\mathrm{P},X)$, all of them are all identical up to continuous deformation. More precisely, consider the homeomorphisms $h:\mathbb{R}^3\rightarrow\mathbb{R}^3$ such that $\pi\circ{h}$ coincides with $\pi$, and, for any two points $x$ and $y$ from $\mathbb{R}^3$, if $x$ is below $y$ then $h(x)$ is below $h(y)$. In other words, these homeomorphisms fix the first two coordinates and preserve the order of the third coordinate. We consider the blow-up triangulations of $\mathrm{P}$ up to the homeomorphisms among these that further fix boundary of $\mathrm{P}$. Therefore, we will call $K$ the blow-up triangulation of $\mathrm{P}$ \emph{associated} with $(T_i)_{0\leq{i}\leq{k}}$ in the sequel. We will often need to consider the family $(\Pi_i)_{0\leq{i}\leq{k}}$ of topological disks that we used to build $K$, or the corresponding map $\psi$ from $\{1, \ldots, k\}$ to the set of the tetrahedra contained in $K$. In order to do that we will use the following statement, an immediate consequence of the above construction.

\begin{prop}\label{PW.sec.1.lem.1}
Consider a path $(T_i)_{0\leq{i}\leq{k}}$ in $\mathcal{F}(\mathrm{P},X)$ and the associated blow-up triangulation $K$ of $\mathrm{P}$. There exists a bijection $\psi$ from $\{1,\ldots,k\}$ to the set of the tetrahedra in $K$ and a family $(\Pi_i)_{0\leq{i}\leq{k}}$ of topological disks embedded within $\mathbb{R}^3$ in such a way that
\begin{itemize}
\item[(i)] $\pi$ induces homeomorphisms $\Pi_i\rightarrow\mathrm{P}$,
\item[(ii)] $\psi(i)$ is below $\Pi_j$ if and only if $i\leq{j}$,
\item[(iii)] if $\sigma$ is a face of $T_i$, then $\pi^{-1}(\sigma)\cap\Pi_i$ is a face of $K$, and
\item[(iv)] the vertices of $\pi\circ\psi(i)$ are exactly the endpoints of the two arcs exchanged by the flip between $T_{i-1}$ and $T_i$.
\end{itemize}
\smallskip
\end{prop}

While a unique blow-up triangulation is associated to a given path in $\mathcal{F}(\mathrm{P},X)$, the inverse is not true. For instance, the two paths shown in Fig.~\ref{PW.sec.1.fig.2} are associated to the same blow-up triangulation $K$ of a convex hexagon, sketched in the center of the figure. This is due to the fact that flips that exchange the diagonals of quadrilaterals with pairwise disjoint interiors can be done in any order in a triangulation, resulting in different paths. For this reason, while $K$ is \emph{the} blow-up triangulation of $\mathrm{P}$ associated to a path $(T_i)_{0\leq{i}\leq{k}}$, we will say that $(T_i)_{0\leq{i}\leq{k}}$ is \emph{a} path associated to $K$. However, all of these paths are between the same pair of triangulations. This pair of triangulations can be identified as follows.

\begin{defn}
We denote by $K^-$ the subset of the faces of a blow-up triangulation $K$ that are not above any other face of $K$ and by $K^+$ the subset of its faces that are not below any other face of $K$.
\end{defn}

By a slight abuse of notation, we denote by $\pi(K^-)$ and $\pi(K^+)$ the triangulations of $\mathrm{P}$ obtained by projecting the arcs of $K^-$ or $K^+$ (but not the triangles or vertices as formally, a triangulation of $\mathrm{P}$ is a set of arcs). 
%Finally, observe that a triangle $\tau$ in $K$ is a face of at most two of the tetrahedra of $K$. In fact, $\tau$ is not incident to a tetrahedron if and only if $\tau$ belongs to both $K^-$ and $K^+$. It is contained in a unique tetrahedron if and only if $\tau$ is contained in exactly one of $K^-$ and $K^+$. Otherwise, $\tau$ is a face of two tetrahedra in $K$.

\section{Properties of blow-up triangulations}\label{PW.sec.1.1}

In this section, we establish some of the properties of blow-up triangulations. We recall that $\mathrm{P}$ is a convex polygon and $X$ is a superset of the vertex set of $\mathrm{P}$ that may further contain flat vertices of $\mathrm{P}$. Throughout the section, $K$ will denote a blow-up triangulation of $\mathrm{P}$ with vertex set $X$. First consider three arcs contained in $K$ and assume that the union of these arcs is a (topological) circle $C$. If an arc of $K$ penetrates $C$ in the following sense, then $C$ cannot be the boundary of a triangle of $K$.

\begin{defn}\label{PW.sec.1.1.defn.1}
Consider a topological circle $C$ obtained as the union of three arcs of a blow-up triangulation $K$. An arc of $K$ \emph{penetrates} $C$ when it is above one of the arcs of $K$ contained in $C$ and below another such arc.
\end{defn}

In order to prove the following lemma, the usual topological notions of the star and link of a face of a complex will be useful. More precisely, consider a face $\sigma$ of a blow-up triangulation $K$. The \emph{star} of $\sigma$ in $K$ is the set $S_K(\sigma)$ of all the faces of $K$ that admit $\sigma$ as one of their faces and the \emph{link} of $F$ in $K$ is the set $L_K(\sigma)$ made of all the faces $\tau$ of $K$ disjoint from $\sigma$ such that $\tau$ is a face of some face of $K$ contained in $S_K(\sigma)$.

\begin{lem}\label{PW.sec.1.1.lem.2}
Consider a topological circle $C$ formed as the union of three arcs of a blow-up triangulation $K$ and let $\beta$ be one of these arcs. If $S_K(\beta)$ does not contain a triangle $\tau$ such that $\pi(C)$ is the boundary of $\pi(\tau)$, then there exists an arc in $L_K(\beta)$ that penetrates $C$.
\end{lem}
\begin{proof}
Denote by $\alpha$ and $\gamma$ the two arcs of $K$ other than $\beta$ that are contained in $C$. Assume that there is no triangle in $S_K(\beta)$ whose image by $\pi$ admits $\pi(C)$ for its boundary. Observe that, under this assumption, the vertex $b$ shared by $\alpha$ and $\gamma$ cannot be contained in the link of $\beta$ in $K$.

If $\beta$ does not belong to $K^-$, then there is a unique tetrahedron $\sigma^-$ incident to and below $\beta$ as sketched in Fig.~\ref{PW.sec.1.1.fig.1.5}. If the edge $\varepsilon^-$ of $\sigma^-$ that is disjoint from $\beta$ is above $\alpha$ or above $\gamma$, then that edge penetrates $C$ and the lemma is proven (the figure only shows the case when $\varepsilon^-$ is above or below $\gamma$ but the case when it is above or below $\alpha$ is symmetric). Therefore, we assume that $\varepsilon^-$ is not above $\alpha$ and not above $\gamma$. As $\varepsilon^-$ is not incident to $b$, it must be below $\alpha$ or $\gamma$, and the same must be true of exactly one of the facets of $\sigma^-$ incident to $\beta$ (which is the striped triangle in Fig.~\ref{PW.sec.1.1.fig.1.5}). Let us denote by $\tau^-$ that facet. In the case when $\beta$ is contained in $K^-$, we shall just denote by $\tau^-$ the triangle in $K^-$ that is incident to $\beta$ and that lies below $\alpha$ or below $\gamma$. Similarly, we can either find an arc in $L_K(\beta)$ that penetrates $C$ (in which case the lemma is proven) or a triangle $\tau^+$ in $S_K(\beta)$ that lies above $\alpha$ or above $\gamma$. Assuming the latter, we now have two triangles $\tau^-$ and $\tau^+$ in $S_K(\beta)$ such that $\tau^-$ is below $\alpha$ or $\gamma$ and $\tau^+$ is above one of these arcs. We will now show that some arc in $L_K(\beta)$ penetrates $C$.

\begin{figure}
\begin{centering}
\includegraphics[scale=1]{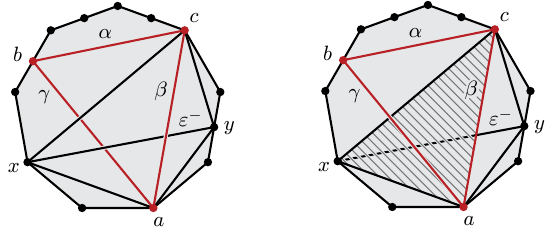}
\caption{The tetrahedron $\sigma^-$ from the proof of Lemma~\ref{PW.sec.1.1.lem.2} (with vertices $a$, $c$, $x$, and $y$) when $\varepsilon^-$ is above $\gamma$ (left) or below $\gamma$ (right). The triangle $\tau^-$ is striped.}\label{PW.sec.1.1.fig.1.5}
\end{centering}
\end{figure}

Consider the triangle in $S_K(\beta)$ below $\alpha$ or $\gamma$ that is above every other such triangle. Similarly, consider the triangle in $S_K(\beta)$ above $\alpha$ or $\gamma$ that is below every other such triangle. There exists a tetrahedron $\sigma$ in $S_K(\beta)$ that is incident to both of these triangles. This follows from the observation that, when two triangles in $S_K(\beta)$ are below or above $\alpha$ or $\gamma$, these triangles must be below or above one another. The tetrahedron $\sigma$ has an edge that belongs to $L_K(\beta)$ and that edge necessarily penetrates $C$, as desired.
\end{proof}

\begin{figure}[b]
\begin{centering}
\includegraphics[scale=1]{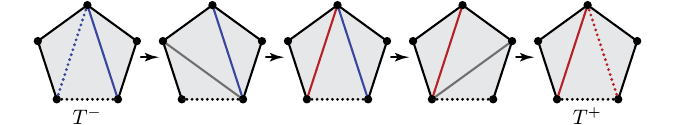}
\caption{A sequence of flips between two triangulations $T^-$ and $T^+$ of a convex pentagon. The dotted arcs, thought of as arcs of the blow-up triangulation associated to that sequence of flips form a topological circle.}\label{PW.sec.1.1.fig.1}
\end{centering}
\end{figure}

\begin{rem}\label{PW.sec.1.1.rem.1}
Not only does Lemma \ref{PW.sec.1.1.lem.2} assume that $C$ is not the boundary of any triangle of $K$ but that $\pi(C)$ is not the boundary of the image by $\pi$ of any of triangle of $K$. This stronger assumption is required. Indeed, consider the sequence of flips depicted in Fig. \ref{PW.sec.1.1.fig.1} between two triangulations $T^-$ and $T^+$ of a convex pentagon $\mathrm{P}$ and the blow-up triangulation $K$ of $\mathrm{P}$ associated to it. The three arcs in $K$ that correspond to the dotted arcs in Fig. \ref{PW.sec.1.1.fig.1} form a topological circle $C$. However, $C$ is not the boundary of a triangle of $K$ and no arc in $K$ penetrates $C$ either. Lemma \ref{PW.sec.1.1.lem.2} fails here because $\pi(C)$ is the boundary of the image by $\pi$ of (three) triangles of $K$.
\end{rem}

\begin{figure}
\begin{centering}
\includegraphics[scale=1]{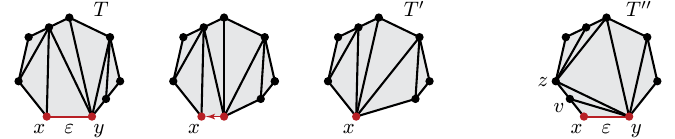}
\caption{A triangulation $T$ (left) and the triangulation $T'$ obtained by contracting $\varepsilon$ to $x$ (center). The edge $\varepsilon$ cannot be contracted to $x$ in the triangulation $T''$ (right).}\label{PW.sec.1.1.fig.2}
\end{centering}
\end{figure}

Let us now translate in terms of blow-up triangulations some of the tools introduced in~\cite{Pournin2014} to estimate distances in flip-graphs. In order to do that, we first need to define two operations. The first operation, already used in~\cite{Pournin2014} consists in contracting an arc between consecutive points of $X$ within a triangulation (recall that $x$ and $y$ are consecutive in $X$ when the line segment between them is contained in the boundary of $\mathrm{P}$ and does not contain any point of $X$ other than $x$ and $y$). The second operation is related: it consists in pulling a tetrahedron, a triangle, or an arc to a point from $X$ and may be carried out within a blow-up triangulation of $\mathrm{P}$ in order to build another blow-up triangulation of the same polygon. We begin by describing edge contractions and illustrate them on the left of Fig. \ref{PW.sec.1.1.fig.2}.

\begin{defn}
Consider a triangulation $T$ of $\mathrm{P}$ with vertex set $X$ and an arc $\varepsilon$ between two consecutive points $x$ and $y$ of $X$. If the point of $X$ consecutive to $x$ but distinct from $y$ is not a flat vertex of $\mathrm{P}$, then we can build a triangulation $T'$ of the convex hull of $X\mathord{\setminus}\{y\}$ by
 \begin{enumerate}
 \item[(i)] removing from $T$ the arc with vertices $x$ and $y$ and
 \item[(ii)] replacing every remaining arc from $T$ between $y$ and a vertex $z$ by as arc between $x$ and a vertex $z$.
 \end{enumerate}
 
We say that $T'$ is obtained from $T$ by \emph{contracting} $\varepsilon$ to $x$.
\end{defn}

Note that, if the point $v$ of $X$ consecutive to $x$ and distinct from $y$ is a flat vertex of $\mathrm{P}$, one cannot always built $T'$ from $T$, as shown on the right of Fig. \ref{PW.sec.1.1.fig.2}. Indeed, in this case, one of the arcs of $T$ may be sent by the contraction to an arc that contains $v$ in its interior, as for instance the arc with vertices $y$ and $z$ on the right of Fig. \ref{PW.sec.1.1.fig.2}. 

We now explain how the pulling operation affects a tetrahedron.

\begin{defn}\label{PW.sec.1.1.defn.2}
Consider a tetrahedron $\sigma$ from a blow-up triangulation $K$ and a tetrahedron $\sigma'$ from a blow-up triangulation $K'$. We say that $\sigma'$ is obtained by \emph{pulling} $\sigma$ to a vertex $x$ of $K'$ when $\sigma$ is incident to a triangle $\tau$ of $K$ and $\sigma'$ to a triangle $\tau'$ of $K'$ such that
\begin{itemize}
\item[(i)] $x$ is not a vertex of $\tau'$ or $\sigma$,
\item[(ii)] the images of $\tau$ and $\tau'$ by $\pi$ coincide, and
\item[(iii)] $\sigma'$ is below $\tau'$ if and only if $\sigma$ is below $\tau$.
\end{itemize}
\end{defn}

\begin{figure}[b]
\begin{centering}
\includegraphics[scale=1]{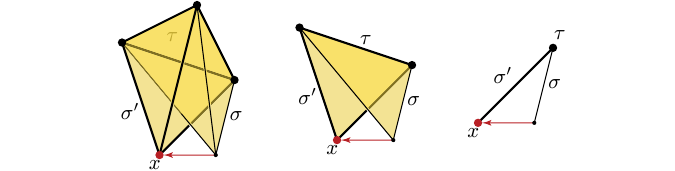}
\caption{The operation of pulling a tetrahedron (left), a triangle (center), and an arc (right) to a point $x$.}\label{PW.sec.1.1.fig.3}
\end{centering}
\end{figure}

In the case of an arc or a triangle, the pulling operation is similar but simpler to describe. Consider two such objects $\sigma$ and $\sigma'$ of the same dimension. We say that $\sigma'$ is obtained by pulling $\sigma$ to $x$ when $x$ is a vertex of $\sigma'$ but not a vertex of $\sigma$ and a facet $\tau$ of $\sigma$ has the same image by $\pi$ than the facet $\tau'$ of $\sigma'$ that does not admit $x$ as a vertex. This operation is illustrated at the center of Fig. \ref{PW.sec.1.1.fig.3} when $\sigma$ is a triangle and on the left of the figure when $\sigma$ is an arc. One may modify a blow-up triangulation $K$ with vertex set $X$ into another blow-up triangulation by removing a subset of its faces and then by pulling to a point $x$ of $X$ a subset of the remaining faces of $K$. When a face of $K$ is pulled to $x$, its lower-dimensional faces that are not contained in $\tau$ are pulled to $x$ along with the face itself. In addition, the faces of $K$ that are not removed or pulled to $x$ remain unaffected and belong to the resulting blow-up triangulation. %We will illustrate both operations in the following of this section in Section~\ref{PW.sec.1.5}.
The proof of the following lemma explains how the edge contraction and pulling operations are related.

\begin{lem}\label{PW.sec.1.1.lem.4}
Consider an arc $\varepsilon$ between two consecutive points $x$ and $y$ of $X$. Let $\mathrm{P}'$ be the convex hull of $X\mathord{\setminus}\{y\}$. If the point of $X$ consecutive to $x$ and distinct from $y$ is not flat, then a blow-up triangulation of $\mathrm{P}'$ with vertex set $X\mathord{\setminus}\{y\}$ can be obtained $K$ by removing all the faces incident to $\varepsilon$, and by pulling to $x$ the remaining faces that are incident to $y$.
\end{lem}
\begin{proof}
Consider a path $(T_i)_{0\leq{i}\leq{k}}$ in $\mathcal{F}(\mathrm{P},X)$ associated to $K$ as well as the corresponding bijection $\psi$ from $\{1, \ldots, k\}$ to the set of the tetrahedra in $K$ provided by Proposition \ref{PW.sec.1.lem.1}. Under the assumption that the point of $X$ consecutive to $x$ but distinct from $y$ is not a flat vertex, we can consider the triangulation $T'_i$ of $\mathrm{P}'$ obtained by contracting $\varepsilon$ to $x$ within $T_i$. Observe that $T'_{i-1}$ and $T'_i$ are identical when the quadrilateral $\pi\circ\psi(i)$ admits $\varepsilon$ as an edge. Moreover, when $\pi\circ\psi(i)$ does not admits $\varepsilon$ as an edge, then $T'_{i-1}$ and $T'_i$ are still related by a flip. Hence, one obtains a path in $\mathcal{F}(\mathrm{P}',X\mathord{\setminus}\{y\})$ by removing consecutive duplicates from $(T'_i)_{0\leq{i}\leq{k}}$. By construction, the blow-up triangulation of $\mathrm{P}'$ with vertex set $X\mathord{\setminus}\{y\}$ associated to this path is obtained by removing from $K$ the faces that are incident to $\varepsilon$ and by pulling to $x$ all the remaining faces incident to $y$ as desired.
\end{proof}

The following lemma allows to bound the number of tetrahedra in $K$ that are incident to an arc between two consecutive points of $X$. It corresponds to Theorem~5 from \cite{Pournin2014} but in terms of blow-up triangulations.

\begin{lem}\label{PW.sec.1.1.lem.3}
Consider a point $z$ in $X$ and the two arcs $\varepsilon$ and $\varepsilon'$ between $z$ and the points in $X$ consecutive to $z$. If $z$ is incident to at least two arcs from $K^-\mathord{\setminus}K^+$ and is not incident to any of the arcs from $K^+\mathord{\setminus}K^-$, then $\varepsilon$ and $\varepsilon'$ cannot both be incident to only one tetrahedron of $K$.
\end{lem}
\begin{proof}
Assume that $z$ is not incident to any arc from $K^+\mathord{\setminus}K^-$. In that case $K^+$ contains a triangle $\tau$ that admits $\varepsilon$ and $\varepsilon'$ as two of its edges. Let us assume that only one tetrahedron $\sigma_\varepsilon$ in $K$ is incident to $\varepsilon$. In this case, $\sigma_\varepsilon$ must admit $\tau$ as a facet. Therefore, $\sigma_\varepsilon$ is also incident to $\varepsilon'$. Denote by $a$ the vertex of $\sigma_x$ that is not a vertex of $\tau$. As $\sigma_\varepsilon$ is the only tetrahedron in $K$ incident to $\varepsilon$, its facet incident to $\varepsilon$ and to $a$ is contained in $K^-$. Now assume that $z$ is incident to at least two arcs from $K^-\mathord{\setminus}K^+$. In this case the triangle in $K^-$ incident to $\varepsilon'$ cannot admit $a$ as a vertex or $\varepsilon$ as an edge. Hence, the tetrahedron $\sigma_{\varepsilon'}$ in $K$ incident to that triangle is distinct from $\sigma_\varepsilon$. Since $\sigma_\varepsilon$ and $\sigma_{\varepsilon'}$ are both incident to $\varepsilon'$, this proves the lemma.
\end{proof}

Lemma \ref{PW.sec.1.1.lem.3} makes it possible to lower bound the number of tetrahedra in a blow-up triangulation $K$ in terms of the number of tetrahedra in a blow-up triangulation obtained from $K$ by the pulling operation described by Lemma~\ref{PW.sec.1.1.lem.4}. In Sections~\ref{PW.sec.3} and~\ref{PW.sec.4}, we will need a consequence of these lemmas that deals with sequences of pulling operations. Let us fix an arc $\delta$ of $K^+$ and denote its two vertices by $x$ and $y$. If $x$ and $y$ are not consecutive points of $X$, then we denote by $\Xi$ and $\Omega$ the two convex polygons obtained by cutting $\mathrm{P}$ along $\pi(\delta)$. If however $x$ and $y$ are consecutive points of $X$, then $\delta$ is contained in the boundary of $\mathrm{P}$ and we take $\mathrm{P}$ for $\Omega$ and $\delta$ for $\Xi$. %The set of the vertices of $\Sigma$ will be denoted by $\mathcal{V}$ and the set of the vertices of $\Sigma^\delta$ by $\mathcal{V}^\delta$.
In order to prove the announced consequence of Lemmas \ref{PW.sec.1.1.lem.4} and \ref{PW.sec.1.1.lem.3}, we will need the following straightforward property.

\begin{prop}\label{PW.sec.1.1.prop.1}
If $\Xi$ has at least three vertices, then it admits a vertex that is not incident to any arc from $K^+\mathord{\setminus}K^-$.
\end{prop}
\begin{proof}
Assume that $\Xi$ has at least three vertices and consider the triangulation $U$ of $\Xi$ contained in $\pi(K^+)$. If $U$ contains a single triangle, then consider the vertex $z$ of that triangle other than $x$ and $y$. In that case, the only two arcs in $K^+$ that are incident to $z$ are the two edges of $\Xi$ incident to $z$. Therefore, $z$ cannot be incident to any arc from $K^+\mathord{\setminus}K^-$.

Now assume that $U$ contains more than just one triangle, and let us recall the well-known property that in this case, at least two triangles of $U$ have two of their edges in the boundary of $\Xi$ (such triangles are often called \emph{ears} \cite{ParlierPournin2017,Pournin2019}). At least one of these triangles does not admit $\pi(\delta)$ as an edge. The vertex shared by the two edges of that triangle that are contained in the boundary of $\Xi$ cannot be incident to any arc from $K^+\mathord{\setminus}K^-$. 
\end{proof}

\begin{figure}
\begin{centering}
\includegraphics[scale=1]{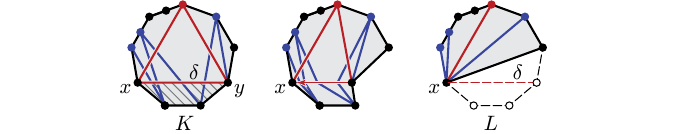}
\caption{Partial sketch of a blow-up triangulation $K$ that satisfies the conditions of Theorem \ref{PW.sec.1.1.thm.1} (left) and of the blow-up triangulation $L$ provided by the theorem (right).}\label{PW.sec.1.1.fig.4}
\end{centering}
\end{figure}

We now state and prove the announced result, that can be thought of as a generalization to blow-up triangulations of Corollary 2 from \cite{Pournin2014}. Under appropriate conditions on $K$, it provides a blow-up triangulation $L$ of the convex hull $\Lambda$ of $(X\mathord{\setminus}\{y\})\cap\Omega$ whose number of tetrahedra is less than the number of tetrahedra of $K$ by a controlled amount. A blow-up triangulation $K$ that satisfies the required conditions is sketched on the left of Fig. \ref{PW.sec.1.1.fig.4}, and the resulting blow-up triangulation $L$ is sketched on the right of that figure (where the portion of $\mathrm{P}$ outside of $\Lambda$ and the arc $\delta$ are drawn using dashed lines). On the left of the figure, the polygon $\Xi$ is striped. The arcs in $K^-\mathord{\setminus}K^+$ incident to a point of $X$ in $\Xi$ are colored blue as well as the links in $K^-$ of the edges of $\mathrm{P}$ contained in $\Xi$. Note that each of these links is made of a single vertex. The edges of the triangle in $K^+\mathord{\setminus}K^-$ incident to $\delta$ and whose image by $\pi$ is contained in $\Omega$ are colored red as well as its vertex opposite $\delta$. Note that this vertex is the portion of the link of $\delta$ in $K^+$ that is disjoint from $\Xi$. %Note that, when $\Xi$ is a polygon (and not a line segment), the link of $\delta$ in $K^+$ is made of two points, one in $\Xi$ and one in $\Omega$, but we only show the point contained in $\Omega$ in the figure because the required condition is always satisfied for the other point.
The way the faces of $K$ will be pulled to $x$ in order to build $L$ is sketched at the center of Fig. \ref{PW.sec.1.1.fig.4}.

\begin{thm}\label{PW.sec.1.1.thm.1}
Assume that the convex hull $\Lambda$ of $(X\mathord{\setminus}\{y\})\cap\Omega$ is not a line segment. Further assume that the flat vertices of $\mathrm{P}$ contained in $X$ all belong to $\Omega$ and are not consecutive to $x$ or $y$. If the links in $K^-$ of the edges of $\mathrm{P}$ contained in $\Xi$ are pairwise distinct and all of them are disjoint from both $\Xi$ and the link of $\delta$ in $K^+$, then there exists a blow-up triangulation $L$ of $\Lambda$ with vertex set $(X\mathord{\setminus}\{y\})\cap\Omega$ such that
\begin{itemize}
\item[(i)] the number of tetrahedra contained in $L$ is less than the number of tetrahedra contained in $K$ by at least $2|X\cap\Xi|-3$ and
\item[(ii)] $L$ is obtained from $K$ by removing any face incident to more than one vertex of $\Xi$ and by pulling to $x$ all of the remaining faces that are incident to a vertex of $\Xi$ other than $x$.
\end{itemize}
\end{thm}
\begin{proof}
The proof is by induction on the number of points in $X\mathord{\setminus}\Omega$. If there is no such point, then $\delta$ is an edge of $\mathrm{P}$ and we apply Lemma~\ref{PW.sec.1.1.lem.4}. According to this lemma, removing from $K$ all the faces incident to $\delta$ and pulling to $x$ all the remaining faces that are incident to $y$ results in a blow-up triangulation $L$ of $\Lambda$ with vertex set $X\mathord{\setminus}\{y\}$. The disjointness condition between the links in $K^-$ of the edges of $\mathrm{P}$ contained in $\Xi$ and the link of $\delta$ in $K^+$ amounts to ask that the third vertex of the triangles incident to $\delta$ in $K^-$ and in $K^+$ are distinct. It follows that $\delta$ is incident to at least two distinct triangles of $K$, and therefore to at least one of its tetrahedra. As a consequence, $L$ contains at least one tetrahedron less than $K$, as desired.

Now assume that $X\mathord{\setminus}\Omega$ is non-empty or equivalently that $\Xi$ has at least three vertices. By Proposition~\ref{PW.sec.1.1.prop.1}, some vertex $c$ of $\Xi$ is not incident to any arc from $K^+\mathord{\setminus}K^-$. As $\delta$ belongs to $K^+\mathord{\setminus}K^-$ and admits $x$ and $y$ as its vertices, $c$ is distinct from $x$ and $y$. Hence, the two edges $\varepsilon$ and $\varepsilon'$ of $\mathrm{P}$ incident to $c$ are contained in $\Xi$. Under the assumption that the links in $K^-$ of $\varepsilon$ and $\varepsilon'$ are distinct, $c$ is incident to at least two distinct arcs from $K^-\mathord{\setminus}K^+$. Therefore, by Lemma~\ref{PW.sec.1.1.lem.3}, $\varepsilon$ and $\varepsilon'$ cannot both be incident to only one tetrahedron of $K$. We assume without loss of generality that $\varepsilon$ is an edge of at least two tetrahedra of $K$. Denote by $\mathrm{P}'$ the convex hull of $X\mathord{\setminus}\{c\}$ and by $a$ the vertex of $\varepsilon$ distinct from $c$. According to Lemma~\ref{PW.sec.1.1.lem.4}, a blow-up triangulation $M$ of $\mathrm{P}'$ with vertex set $X\mathord{\setminus}\{c\}$ can be obtained from $K$ by removing all the faces incident to $\varepsilon$ and by pulling to $a$ the remaining faces that are incident to $c$. Since $K$ contains at least two tetrahedra that admit $\varepsilon$ as an edge, the number of tetrahedra contained in $M$ is less by at least two than the number of tetrahedra contained in $K$.

Observe that $M$ satisfies the conditions of the lemma with respect to $\mathrm{P}'$ as soon as $K$ does with respect to $\mathrm{P}$. Note, in particular that the links in $M^-$ of the edges of $\mathrm{P}'$ contained in $\Xi$ form a subset of the links in $K^-$ of the edges of $\mathrm{P}$ contained in $\Xi$. Moreover, the number of vertices of $\mathrm{P}'$ that do not belong to $\Omega$ is less by one than the number of points in $X\mathord{\setminus}\Omega$. Hence, by induction, there exists a blow-up triangulation $L$ of $\Lambda$ whose number of tetrahedra is less than the number of tetrahedra in $M$ by at least $2|X\cap\Xi|-5$. It is obtained from $M$ by removing any face incident to more than one vertex of $\Xi$, and by pulling to $x$ the remaining faces that are incident to a vertex of $\Xi$ other than $x$. It immediately follows that the number of tetrahedra in $L$ is less than the number of tetrahedra in $K$ by at least $2|X\cap\Xi|-3$. Moreover, by construction, $L$ is obtained from $K$ by removing any face incident to more than one vertex of $\Xi$ and by pulling to $x$ all of the remaining faces that are incident to a vertex of $\Xi$ other than $x$, as desired.
\end{proof}

\section{Revisiting the projections related to arcs}\label{PW.sec.1.5}
 
\begin{figure}
\begin{centering}
\includegraphics[scale=1]{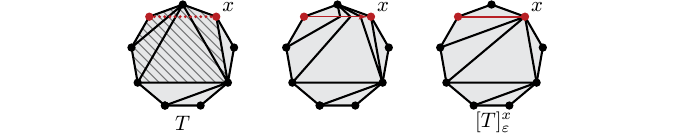}
\caption{The projection of a triangulation $T$ (left) to the triangulation $[T]_\varepsilon^x$ (right), where $\varepsilon$ is the arc colored red. The pulling operation is shown in the center of the figure.}\label{PW.sec.1.5.fig.1}
\end{centering}
\end{figure}

Let us assume for a moment that $X$ does not contain any flat vertex of $\mathrm{P}$ or equivalently, that it is the vertex set of $\mathrm{P}$. We recall how the projection from $\mathcal{F}(\mathrm{P},X)$ to $\mathcal{F}_\varepsilon(\mathrm{P},X)$ introduced in \cite{SleatorTarjanThurston1988} works where $\varepsilon$ is an arc between two points of $X$. Consider a triangulation $T$ of $\mathrm{P}$ with vertex set $X$. If $T$ admits $\varepsilon$ as an arc, then $T$ already belongs to $\mathcal{F}_\varepsilon(\mathrm{P},X)$ and coincides with its projection, that we shall denote by $[T]_\varepsilon^x$. If $\varepsilon$ is not contained in $T$, then the union of the triangles in $T$ whose interior is non-disjoint from $\varepsilon$ is a convex polygon $\mathrm{P}'$. In other words, $T$ admits as a subset a triangulation $T'$ of $\mathrm{P}'$. The polygon $\mathrm{P}'$ is striped on the left of Fig. \ref{PW.sec.1.5.fig.1}.

Now recall that the \emph{comb} triangulation of a convex polygon at one of its vertices $x$ is the triangulation that contains arcs incident to $x$ and to each of the other vertices of that polygon. One can build a new triangulation $[T]_\varepsilon^x$ of $\mathrm{P}$ by replacing within $T$ the arcs contained in $T'$ by the comb triangulation of $\mathrm{P}'$ at one of the vertices $x$ of $\varepsilon$, as illustrated in Fig. \ref{PW.sec.1.5.fig.1}. The triangulation $[T]_\varepsilon^x$ is the projection of $T$ to $\mathcal{F}_\varepsilon(\mathrm{P},X)$ described in \cite{SleatorTarjanThurston1988}.

It will be useful to remember that there are two different possible projections to $\mathcal{F}_\varepsilon(\mathrm{P},X)$ depending on which vertex $x$ of $\varepsilon$ is chosen for the comb triangulations. An important property of these projections is that two triangulations related by a flip will be projected either to the same triangulation, or to two triangulations that are still related by a flip \cite{SleatorTarjanThurston1988}. In fact the projections of two triangulations related by a flip coincide if and only if that flip exchanges two arcs that cross $\varepsilon$ or an arc that crosses $\varepsilon$ and an arc incident to $x$. The following is a consequence of these remarks.

\begin{lem}[Lemma 3(b) in \cite{SleatorTarjanThurston1988}]\label{PW.sec.1.5.lem.1}
If $X$ does not contain any flat vertex, then $\mathcal{F}_\varepsilon(\mathrm{P},X)$ is a strongly convex subgraph of $\mathcal{F}(\mathrm{P},X)$. 
\end{lem}

Informally, this projection to $\mathcal{F}_\varepsilon(\mathrm{P},X)$ can be thought of as pulling to $x$ the arcs that cross $\varepsilon$ in a triangulation---as if these arcs were rubber bands---as shown at the center of Fig. \ref{PW.sec.1.5.fig.1} and then inserting the arc $\varepsilon$ as a diagonal of the quadrilateral created by this pulling operation. If $K$ is the blow-up triangulation associated to a path in $\mathcal{F}(\mathrm{P},X)$ and $\varepsilon$ a line segment between two points of $X$, the projection of that path to $\mathcal{F}_{\pi(\varepsilon)}(\mathrm{P},X)$ can also be thought of as an operation within $K$. That operation consists in removing a subset of the tetrahedra and pulling to $x$ another subset in the sense of Definition~\ref{PW.sec.1.1.defn.2}. The tetrahedra that are not removed or pulled to $x$ are unaffected and belong to the resulting blow-up triangulation.

\begin{thm}\label{PW.sec.1.5.thm.1}
Assume that $X$ is the vertex set of $\mathrm{P}$ and consider a line segment $\varepsilon$ between two points of $X$. If $x$ is one of the vertices of $\varepsilon$, then there exists some blow-up triangulation $N$ of $\mathrm{P}$ such that
\begin{itemize}
\item[(i)] $\pi(N^-)$ is equal to $[\pi(K^-)]_\varepsilon^x$ and $\pi(N^+)$ to $[\pi(K^+)]_\varepsilon^x$,
\item[(ii)] $\pi^{-1}(\varepsilon)\cap{N^-}$ coincides with $\pi^{-1}(\varepsilon)\cap{N^+}$, and
\item[(iii)] $N$ is obtained from $K$ by removing the tetrahedra $\sigma$ such that the diagonals of $\pi(\sigma)$ both cross $\varepsilon$ or one of them crosses $\varepsilon$ and the other is incident to $x$ and by pulling to $x$ all the remaining tetrahedra $\sigma$ of $K$ such that a diagonal of $\pi(\sigma)$ crosses $\varepsilon$.
\end{itemize}
\smallskip
\end{thm}
\begin{proof}
Consider a path $(T_i)_{0\leq{i}\leq{k}}$ in $\mathcal{F}(\mathrm{P},X)$ associated to $K$ and denote by $\psi$ the bijection from $\{1,...,k\}$ to the set of the tetrahedra in $K$ that corresponds to that path via Proposition \ref{PW.sec.1.lem.1}. 

\begin{figure}[b]
\begin{centering}
\includegraphics[scale=1]{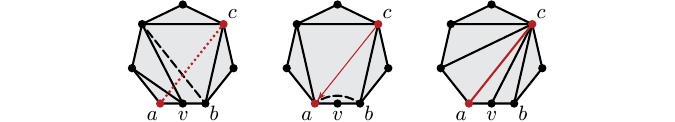}
\caption{A triangulation $T$ (left) and its projection $[T]_\varepsilon^c$ (right), where $\varepsilon$ is colored red. Pulling the dashed arc to $a$ results in an arc that contains $v$ in its interior (center).}\label{PW.sec.1.5.fig.2}
\end{centering}
\end{figure}

% and $(\Pi_i)_{0\leq{i}\leq{k}}$ the family of disks associated to that path by the lemma.
As we have mentioned, the triangulations $[T_{i-1}]_\varepsilon^x$ and $[T_i]_\varepsilon^x$ are either identical or related by a flip. Consider the path from $[\pi(K^-)]_\varepsilon^x$ to $[\pi(K^+)]_\varepsilon^x$  in $\mathcal{F}_\varepsilon(\mathrm{P},X)$ obtained by removing consecutive duplicates from the sequence $([T_i]_\varepsilon^x)_{0\leq{i}\leq{k}}$ as well as the blow-up triangulation $N$ of $\mathrm{P}$ associated to that path. Let us review each of the flips in $(T_i)_{0\leq{i}\leq{k}}$. When a single diagonal of $\pi\circ\psi(i)$ crosses $\varepsilon$ and the other is not incident to $x$, the quadrilateral whose diagonals are exchanged by the flip between $[T_{i-1}]_\varepsilon^x$ and $[T_i]_\varepsilon^x$ shares exactly three vertices with $\pi\circ\psi(i)$ and its fourth vertex is $x$. The tetrahedron in $N$ corresponding to that flip is therefore obtained by pulling $\psi(i)$ to $x$. When $\pi\circ\psi(i)$ does not have a diagonal that crosses $\varepsilon$, the flip that changes $T_{i-1}$ into $T_i$ exchanges the same arcs as the one that changes $[T_{i-1}]_\varepsilon^x$ into $[T_i]_\varepsilon^x$ and the tetrahedron corresponding to this flip is the same in $K$ and $N$. The only remaining case is when the diagonals of $\pi\circ\psi(i)$ both cross $\varepsilon$ or one of them does while the other is incident to $x$. In that case $[T_{i-1}]_\varepsilon^x$ and $[T_i]_\varepsilon^x$ coincide which means that there is no tetrahedron in $N$ that corresponds to $\psi(i)$. Since all the triangulations $[T_i]_\varepsilon^x$ contain $\varepsilon$, $N$ is pinched along its arc whose image by $\pi$ is $\varepsilon$ which completes the proof.
\end{proof}

The strategy in the proof of Theorem \ref{PW.sec.1.5.thm.1} is to consider a path in $\mathcal{F}(\mathrm{P},X)$ associated to a blow-up triangulation, to modify that path, and to consider the blow-up triangulation of $\mathrm{P}$ associated to that modified path. We will use the same general strategy in the next sections. While the properties of the projection described above can be derived from only looking at paths in $\mathcal{F}(\mathrm{P},X)$, the interplay between the combinatorics of these paths and the geometry of the corresponding blow-up triangulations will become instrumental in our proofs of Theorem \ref{PW.sec.0.thm.0} and its variants.

Let us now allow $X$ to contain flat vertices of $\mathrm{P}$. First assume that $X$ contains a single flat vertex $v$ in addition to the vertices of $\mathrm{P}$. Denote by $a$ and $b$ the points of $X$ consecutive to $v$. Pick a Euclidean line segment $\varepsilon$ between two points of $X$ but that does not contain $v$ in its interior. Further consider a vertex $x$ of $\varepsilon$, and a triangulation $T$ of $\mathrm{P}$ with vertex set $X$. A triangulation $[T]_\varepsilon^x$ can be defined just as in the case when $X$ does not contain a flat vertex except if $x$ is equal to $a$ or $b$ and some arc $\delta$ in $T$ that crosses $\varepsilon$ is incident to the point among $a$ and $b$ that is distinct from $x$.

This situation is depicted in Fig. \ref{PW.sec.1.5.fig.2}, in the case when $x$ is equal to $a$ and $\delta$ is the dashed arc. The reason why the projection fails here, is that pulling $\delta$ to $a$ creates an arc that contains $v$ in its interior as shown in the center of the figure. However, since $v$ is the only flat vertex of $\mathrm{P}$, one can always construct the triangulation $[T]_\varepsilon^c$ shown on the right of Fig.~\ref{PW.sec.1.5.fig.2} by using the vertex $c$ of $\varepsilon$ that is not consecutive to $v$ instead of $x$. In other words, there is still a way to project any path in $\mathcal{F}(\mathrm{P},X)$ to $\mathcal{F}_\varepsilon(\mathrm{P},X)$. In this case, the same argument as for Lemma \ref{PW.sec.1.5.lem.1} can be used and we obtain the following.

\begin{lem}\label{PW.sec.1.5.lem.2}
If $X$ is obtained by adding a single flat vertex to the vertex set of $\mathrm{P}$, then $\mathcal{F}_\varepsilon(\mathrm{P},X)$ is a strongly convex subgraph of $\mathcal{F}(\mathrm{P},X)$. 
\end{lem}

\begin{figure}[b]
\begin{centering}
\includegraphics[scale=1]{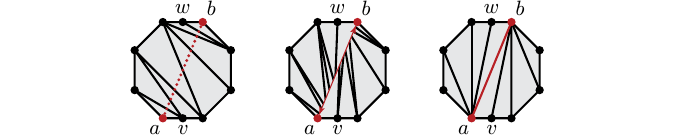}
\caption{A projection (right) of the triangulation $T$ (left) with respect to the arc $\varepsilon$ colored red by pulling (center) the arcs of $T$ crossing $\varepsilon$ to both of its extremities.}\label{PW.sec.1.5.fig.3}
\end{centering}
\end{figure}

Now assume that $X$ is made of the vertices of $\mathrm{P}$ and exactly two flat vertices $v$ and $w$. If $v$ and $w$ are consecutive points of $X$, then the situation is essentially similar to that of Lemma~\ref{PW.sec.1.5.lem.2}. Indeed, in this case, it is always possible to project a path to $\mathcal{F}_\varepsilon(\mathrm{P},X)$ by picking for $x$ a vertex of $\varepsilon$ that is not consecutive to $v$ or $w$. However, if $v$ and $w$ are not consecutive, then $\varepsilon$ can be such that one of its vertices is consecutive to $v$ and the other to $w$. If $w$ is not on the same side of $\varepsilon$ as $v$ then the modified projection illustrated in Fig.~\ref{PW.sec.1.5.fig.3} is still possible. Let us explain how that projection works.

Denote by $a$ the vertex of $\varepsilon$ adjacent to $v$ and by $b$ the vertex of $\varepsilon$ adjacent to $w$. Consider a triangulation $T$ of $\mathrm{P}$ with vertex set $X$ and the convex polygon $\mathrm{P}'$ obtained as the union of the triangles of $T$ whose interior is non-disjoint from $\varepsilon$. Now replace the arcs of $T$ contained in $\mathrm{P}'$ with the triangulation of $\mathrm{P}'$ that contains $\varepsilon$, an arc between $b$ and each of the points in $X\cap\mathrm{P}'$ on the same side of $\varepsilon$ as $v$, and an arc between $a$ and each of the points in $X\cap\mathrm{P}'$ on the same side of $\varepsilon$ than $w$. This projection is illustrated in Fig.~\ref{PW.sec.1.5.fig.3}. Again, the projections of two triangulations that are related by a flip are still either identical or related by a flip and $\mathcal{F}_\varepsilon(\mathrm{P},X)$ remains strongly convex. If, however, $v$ and $w$ are on the same side of $\varepsilon$, then we shall see in Section~\ref{PW.sec.3} that $\mathcal{F}_\varepsilon(\mathrm{P},X)$ is no longer always strongly convex. Finally note that the projection illustrated in Fig. \ref{PW.sec.1.5.fig.3} can still be performed when $X$ contains an arbitrary number of flat vertices of $\mathrm{P}$, provided that $X$ does not have two flat vertices on the same side of $\varepsilon$ each consecutive to a vertex of $\varepsilon$. We obtain the following as a consequence.

\begin{lem}\label{PW.sec.1.5.lem.3}
Consider a set $X$ obtained from the vertex set of $\mathrm{P}$ by adding finitely-many flat vertices of $\mathrm{P}$ and a line segment $\varepsilon$ between two points of $X$ that does not contain any point of $X$ in its interior. If $X$ does not contain two flat vertices on the same side of $\varepsilon$ each consecutive to one of the vertices of $\varepsilon$, then $\mathcal{F}_\varepsilon(\mathrm{P},X)$ is a strongly convex subgraph of $\mathcal{F}(\mathrm{P},X)$.
\end{lem}

It is noteworthy that, as the projection shown in Fig. \ref{PW.sec.1.5.fig.3}, our proofs in the next section will consist in partly re-triangulating a blow-up triangulation of a convex polygon with flat vertices by pulling some of the tetrahedra it contains to one of its vertices and some others to a different vertex.

\section{A projection related to triangles}\label{PW.sec.2}

\begin{figure}[b]
\begin{centering}
\includegraphics[scale=1]{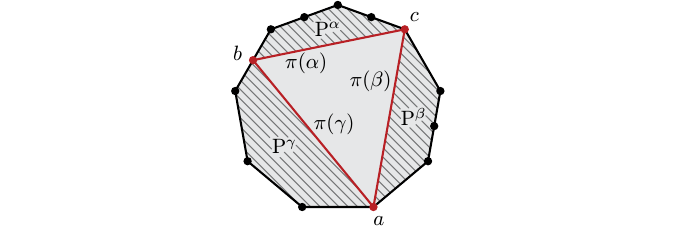}
\caption{The polygons $\mathrm{P}^\alpha$, $\mathrm{P}^\beta$, and $\mathrm{P}^\gamma$.}\label{PW.sec.2.fig.0}
\end{centering}
\end{figure}

In this section, $K$ is a blow-up triangulation of a convex Euclidean polygon $\mathrm{P}$. We will first assume that the vertex set $X$ of $K$ coincides with the vertex set of $\mathrm{P}$ and discuss the case when $X$ further contains flat vertices of $\mathrm{P}$ at the end of the section. Throughout the section, we also fix a topological circle $C$ obtained as the union of three arcs $\alpha$, $\beta$, and $\gamma$ of $K$ such that $\alpha$ is contained in $K^-$ and $\gamma$ in $K^+$. The vertices of these arcs will be denoted by $a$, $b$, and $c$ in such a way that $a$ does not belong to $\alpha$, $b$ does not belong to $\beta$, and $c$ does not belong to $\gamma$. If $\varepsilon$ is one of the arcs $\alpha$, $\beta$, or $\gamma$, we will denote by $\mathrm{P}^\varepsilon$ the portion of $\mathrm{P}$ sketched in Fig.~\ref{PW.sec.2.fig.0}, made of the points of $\mathrm{P}$ that do not lie on the same side of $\pi(\varepsilon)$ as the interior of the triangle with vertices $a$, $b$, and $c$. Note that $\pi(\varepsilon)$ is an edge of $\mathrm{P}^\varepsilon$. Moreover, $\mathrm{P}^\varepsilon$ is a convex polygon when $\pi(\varepsilon)$ is not an edge of $\mathrm{P}$ and it is equal to $\pi(\varepsilon)$ otherwise.

We also consider a path $(T_i)_{0\leq{i}\leq{k}}$ in $\mathcal{F}(\mathrm{P},X)$ associated to $K$. In particular, there exist a bijection $\psi$ from $\{1, ..., k\}$ to the set of the tetrahedra in $K$ and a family $(\Pi_i)_{0\leq{i}\leq{k}}$ of topological disks satisfying assertions (i), (ii), (iii), and (iv) in the statement of Proposition \ref{PW.sec.1.lem.1}. Throughout this section, $q$ is a fixed integer such that $\beta$ is contained in $\Pi_q$. Most of the intermediate results established in this section depend on $q$ and on other objects that we might need to define along the way, but these dependences will disappear at the end of the section when we state the main results. Denote
$$
\left\{
\begin{array}{l}
\Pi^\alpha=\pi^{-1}(\mathrm{P}^\alpha)\cap\Pi_0\mbox{,}\\
\Pi^\beta=\pi^{-1}(\mathrm{P}^\beta)\cap\Pi_q\mbox{, and}\\
\Pi^\gamma=\pi^{-1}(\mathrm{P}^\gamma)\cap\Pi_k\mbox{.}
\end{array}
\right.
$$

In other words, $\Pi^\alpha$, $\Pi^\beta$, and $\Pi^\gamma$ are the portions of $\Pi_0$, $\Pi_q$, and $\Pi_k$ whose orthogonal projections on $\mathbb{R}^2$ are $\mathrm{P}^\alpha$, $\mathrm{P}^\beta$, and $\mathrm{P}^\gamma$. We can define a projection $T\mapsto\langle{T}\rangle_\varepsilon^x$ to $\mathcal{F}_\varepsilon(\mathrm{P},X)$ different from $T\mapsto[T]_\varepsilon^x$ as follows.

\begin{defn}\label{PW.sec.2.defn.0}
Consider a triangulation $T$ of $\mathrm{P}$ with vertex set $X$ and pick an arc $\varepsilon$ among $\alpha$, $\beta$, and $\gamma$. For any vertex $x$ of $\varepsilon$, we set
$$
\langle{T}\rangle_\varepsilon^x=([T]_\varepsilon^x\mathord{\setminus}T')\cup{T''}
$$
where $T'$ is the set of all the arcs of $[T]_\varepsilon^x$ contained in $\mathrm{P}^\varepsilon$ and $T''$ the set of the images by $\pi$ of the arcs of $K$ contained in $\Pi^\varepsilon$.
\end{defn}

\begin{figure}
\begin{centering}
\includegraphics[scale=1]{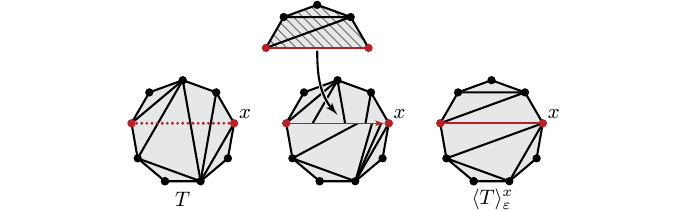}
\caption{A triangulation $T$ of a convex polygon without flat vertices (left) and its projection $\langle{T}\rangle_\varepsilon^x$ (right), where the arc $\pi(\varepsilon)$ is colored red. The striped polygon is $\mathrm{P}^\varepsilon$ (top).}\label{PW.sec.2.fig.1}
\end{centering}
\end{figure}

In other words, $\langle{T}\rangle_\varepsilon^x$ is obtained by first pulling to $x$ the arcs of $T$ that cross $\varepsilon$ as when building $[T]_\varepsilon^x$ from $T$ but then the arcs resulting from this operation and that are contained in $\mathrm{P}^\varepsilon$ are replaced with the images by $\pi$ of the arcs of $K$ contained in $\Pi^\varepsilon$ as sketched in Fig. \ref{PW.sec.2.fig.1}. These images are sketched at the top of the figure where the polygon $\mathrm{P}^\varepsilon$ is striped.

A straightforward but important observation is that, just as for the projections that we discussed in Section~\ref{PW.sec.1.5}, the projection $T\mapsto\langle{T}\rangle_\varepsilon^x$ sends two triangulations related by a flip either to the same triangulation, or to two different triangulations that are still related by a flip. By construction, we obtain the following characterization of these two cases.

\begin{prop}\label{PW.sec.2.prop.1}
Consider two triangulations $T$ and $U$ of $\mathrm{P}$ with vertex set $X$ and assume that they are related by a flip. The triangulations $\langle{T}\rangle_\varepsilon^x$ and $\langle{U}\rangle_\varepsilon^x$ coincide if and only if the arcs exchanged by that flip are both non-disjoint from $\mathrm{P}^\varepsilon\mathord{\setminus}\{y\}$,  where $y$ is the vertex of $\varepsilon$ other than $x$. Otherwise the triangulations $\langle{T}\rangle_\varepsilon^x$ and $\langle{U}\rangle_\varepsilon^x$ are still related by a flip.
\end{prop}

We now consider two subsets of the tetrahedra contained in $K$ that will play an important role in the upcoming proofs.

\begin{defn}\label{PW.sec.2.defn.1}
We call a tetrahedron $\sigma$ contained in $K$ a \emph{lower tetrahedron of $K$} when it admits a facet $\tau$ contained in $K$ such that
\begin{enumerate}
\item[(i)] $\tau$ is not above $\Pi^\beta$,
\item[(ii)] $\tau$ is disjoint from $\mathrm{P}^\alpha\mathord{\setminus}\{b\}$, and
\item[(iii)] if $\tau$ is contained in $\Pi^\beta$, then $\tau$ is above $\sigma$,
\end{enumerate}
and an \emph{upper tetrahedron of $K$} when $\sigma$ has a facet $\tau$ in $K$ satisfying (i), (ii), and (iii) where ``above'' is replaced by ``below'' and $\alpha$ by $\gamma$.
\end{defn}

While the triangle $\tau$ from Definition~\ref{PW.sec.2.defn.1} is not necessarily contained in $\mathbb{R}^2$, its vertices are. As a consequence, the disjointness property in Assertion~(ii) amounts to ask that the only vertex possibly shared by $\tau$ and $\mathrm{P}^\alpha$ is $b$. From now on, we denote by $\mathcal{L}$ the set of the lower tetrahedra of $K$ and by $\mathcal{U}$ the set of its upper tetrahedra. The following lemma explains how the projections $T\mapsto\langle{T}\rangle_\alpha^c$ and $T\mapsto\langle\langle{T}\rangle_\alpha^c\rangle_\beta^c$ are related to $\mathcal{L}$.

\begin{lem}\label{PW.sec.2.lem.3}
The tetrahedron $\psi(i)$ belongs to $\mathcal{L}$ if and only if either
\begin{enumerate}
\item[(i)] $1\leq{i}\leq{q}$ while $\langle{T_{i-1}}\rangle_\alpha^c$ and $\langle{T_i}\rangle_\alpha^c$ are related by a flip or
\item[(ii)] $q<i\leq{k}$ while $\langle\langle{T_{i-1}}\rangle_\alpha^c\rangle_\beta^c$ and $\langle\langle{T_i}\rangle_\alpha^c\rangle_\beta^c$ are related by a flip.
\end{enumerate}
\end{lem}
\begin{proof}
Consider an integer $i$ such that $1\leq{i}\leq{k}$. Let us first assume that $i$ is at most $q$. According to Proposition \ref{PW.sec.2.prop.1}, the triangulations $\langle{T_{i-1}}\rangle_\alpha^c$ and $\langle{T_i}\rangle_\alpha^c$ are related by a flip if and only if at most one vertex of $\psi(i)$ belongs to $\mathrm{P}^\alpha\mathord{\setminus}\{b\}$. Equivalently, $\psi(i)$ has a facet $\tau$ in $K$ that is disjoint from $\mathrm{P}^\alpha\mathord{\setminus}\{b\}$. As $i$ is at most $q$, $\tau$ cannot be above $\Pi^\beta$ and if it is contained in $\Pi^\beta$, then it must be above $\sigma$, which proves the desired statement.

Now assume that $i$ is greater than $q$. Applying Proposition \ref{PW.sec.2.prop.1} twice shows that the triangulations $\langle\langle{T_{i-1}}\rangle_\alpha^c\rangle_\beta^c$ and $\langle\langle{T_i}\rangle_\alpha^c\rangle_\beta^c$ are related by a flip if and only if at least three vertices of $\psi(i)$ are contained in $\mathrm{P}^\gamma$. Equivalently, $\psi(i)$ has a facet $\tau$ in $K$ whose three vertices are in $\mathrm{P}^\gamma$. Since $i$ is greater than $q$, this situation happens if and only if $\tau$ is disjoint from $\mathrm{P}^\alpha\mathord{\setminus}\{b\}$ and not above or contained in $\Pi^\beta$, thus completing the proof.
\end{proof}

In the sequel, we denote by $\mathcal{P}$ the set of the tetrahedra contained in $K$ that are incident to at least one of the arcs of $K$ that penetrate $C$. Note that $\mathcal{P}$ can be empty in the case when $K$ does not have an arc that penetrates $C$. We recall that according to Remark \ref{PW.sec.1.1.rem.1}, this situation can happen even when $C$ is not the boundary of any of the triangles of $K$.

\begin{figure}
\begin{centering}
\includegraphics[scale=1]{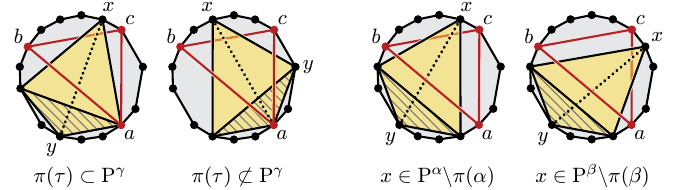}
\caption{The quadrilateral $\pi\circ\psi(i)$ (yellow) when $i$ is at most $q$ (left) and when $i$ is greater than $q$ (right).}\label{PW.sec.2.fig.2}
\end{centering}
\end{figure}

\begin{lem}\label{PW.sec.2.lem.4.1}
Consider an integer $i$ such that $\psi(i)$ is a lower tetrahedron of $K$. That tetrahedron belongs to $\mathcal{P}$ if and only if either
\begin{enumerate}
\item[(i)] $1\leq{i}\leq{q}$ while a diagonal of $\pi\circ\psi(i)$ crosses $\pi(\alpha)$ or
\item[(ii)] $q<i\leq{k}$ and a diagonal of $\pi\circ\psi(i)$ crosses $\pi(\alpha)$ or $\pi(\beta)$.
\end{enumerate}
\end{lem}
\begin{proof}
Let us first treat the case when $i$ is at most $q$. Observe that $\psi(i)$ cannot be above $\gamma$ and, as $i$ is at most $q$, it cannot  be above $\beta$ either. Hence $\psi(i)$ belongs to $\mathcal{P}$ if and only if one of its vertices belongs to $\mathrm{P}^\beta\mathord{\setminus}\pi(\beta)$ or $\mathrm{P}^\gamma\mathord{\setminus}\pi(\gamma)$ and another to $\mathrm{P}^\beta\mathord{\setminus}\pi(\alpha)$. However, by Definition~\ref{PW.sec.2.defn.1}, a facet $\tau$ of $\psi(i)$ in $K$ is disjoint from $\mathrm{P}^\alpha\mathord{\setminus}\{b\}$. Therefore, some vertex $y$ of $\tau$ necessarily belongs to $\mathrm{P}^\beta\mathord{\setminus}\pi(\beta)$ or to $\mathrm{P}^\gamma\mathord{\setminus}\pi(\gamma)$ as illustrated on the left of Fig.~\ref{PW.sec.2.fig.2} depending whether $\pi(\tau)$ is contained in $\mathrm{P}^\gamma$ or not. Hence, it suffices to show that a diagonal of $\pi\circ\psi(i)$ crosses $\pi(\alpha)$ if and only if the vertex $x$ of $\psi(i)$ that is not contained in $\tau$ belongs to $\mathrm{P}^\alpha\mathord{\setminus}\pi(\alpha)$. This equivalence holds because, as $\tau$ is disjoint from $\mathrm{P}^\alpha\mathord{\setminus}\{b\}$, the other vertex of the diagonal of $\pi\circ\psi(i)$ incident to $x$ does not belong to $\mathrm{P}^\alpha$.

Now assume that $i$ is greater than $q$. In this case, $\psi(i)$ cannot be below $\beta$. As $\psi(i)$ cannot be below $\alpha$ either, it is contained in $\mathcal{P}$ if and only if one of its vertices belongs to $\mathrm{P}^\alpha\mathord{\setminus}\pi(\alpha)$ or to $\mathrm{P}^\beta\mathord{\setminus}\pi(\beta)$ and another to $\mathrm{P}^\gamma\mathord{\setminus}\pi(\gamma)$. However since $i$ is greater than $q$, it follows from Definition~\ref{PW.sec.2.defn.1} that some facet $\tau$ or $\psi(i)$ in $K$ must have its three vertices in $\mathrm{P}^\gamma$. In particular, the only vertex of $\psi(i)$ that can belong to $\mathrm{P}^\alpha\mathord{\setminus}\pi(\alpha)$ or to $\mathrm{P}^\beta\mathord{\setminus}\pi(\beta)$ is the one that is not contained in $\tau$, which we denote by $x$. The other vertex $y$ of the diagonal $\varepsilon$ of $\pi\circ\psi(i)$ incident to $x$ is necessarily contained in $\mathrm{P}^\gamma\mathord{\setminus}\pi(\gamma)$ as shown on the right of Fig.~\ref{PW.sec.2.fig.2}. The lemma follows because $\varepsilon$ crosses $\pi(\alpha)$ or $\pi(\beta)$ if and only if $x$ belongs to $\mathrm{P}^\beta\mathord{\setminus}\pi(\alpha)$ or to $\mathrm{P}^\beta\mathord{\setminus}\pi(\beta)$.
\end{proof}

Denote by $T^\circ$ the triangulation of $\mathrm{P}$ whose arcs are the images by $\pi$ of the arcs of $K$ contained in $\Pi^\alpha\cup\Pi^\beta\cup\Pi^\gamma$. Equivalently, $T^\circ$ is the triangulation of $\mathrm{P}$ that shares its arcs contained in $\mathrm{P}^\alpha$ with $T_0$, its arcs contained in $\mathrm{P}^\beta$ with $T_q$, and its arcs contained in $\mathrm{P}^\gamma$ with $T_k$.

\begin{lem}\label{PW.sec.2.lem.5}
There exists a blow-up triangulation $L$ of $\mathrm{P}$ such that
\begin{itemize}
\item[(i)] $\pi(L^-)$ is equal to $\pi(K^-)$ and $\pi(L^+)$ is equal to $T^\circ$, and
\item[(ii)] the tetrahedra of $L$ are obtained from the tetrahedra contained in $\mathcal{L}$ by pulling to $c$ the ones that belong to $\mathcal{P}\cap\mathcal{L}$.
\end{itemize}
\end{lem}
\begin{proof}
Consider an integer $i$ such that $0\leq{i}\leq{k}$. If $i\leq{q}$, denote
$$
T'_i=\langle{T_i}\rangle_\alpha^c\mbox{,}
$$
and if $i>q$, denote
$$
T'_i=\langle\langle{T_i}\rangle_\alpha^c\rangle_\beta^c\mbox{.}
$$

According to Proposition \ref{PW.sec.2.prop.1}, $T'_{i-1}$ and $T'_i$ are either identical or related by a flip. Note that when $i$ is equal to $q+1$, this holds because, as the image by $\pi$ of the arcs of $K$ contained in $\Pi^\beta$ all belong to $T_q$, 
$$
\langle{T_q}\rangle_\alpha^c=\langle\langle{T_q}\rangle_\alpha^c\rangle_\beta^c\mbox{.}
$$

Moreover, $T'_0$ is equal to $\pi(K^-)$ and $T'_k$ to $T^\circ$. As a consequence, one obtains a path in $\mathcal{F}(\mathrm{P},X)$ from $\pi(K^-)$ to $T^\circ$ by removing consecutive duplicates from the sequence $(T'_i)_{0\leq{i}\leq{r}}$. By construction, the set of the tetrahedra of the blow-up triangulation $L$ of $\mathrm{P}$ associated to that path is obtained from $K$ by removing a subset of its tetrahedra, by pulling to $c$ another subset, and by leaving unaffected any other tetrahedron of $K$.

According to Lemma \ref{PW.sec.2.lem.3}, $T'_{i-1}$ and $T'_i$ are related by a flip if and only if $\psi(i)$ belongs to $\mathcal{L}$. Hence, the tetrahedra that are removed when $L$ is obtained from $K$ are exactly the ones that do not belong to $\mathcal{L}$.

Now assume that $\psi(i)$ belongs to $\mathcal{L}$. By Lemma \ref{PW.sec.2.lem.4.1}, the flip between $T_{i-1}$ and $T_i$ exchanges the same two arcs as the flip between $T'_{i-1}$ and $T'_i$ exactly when $\psi(i)$ does not belong to $\mathcal{P}$. As a consequence, $\psi(i)$ belongs to $\mathcal{P}\cap\mathcal{L}$ exactly when the tetrahedron of $L$ that corresponds to the flip between $T'_{i-1}$ and $T'_i$ is obtained by pulling $\psi(i)$ to $c$, as desired.
\end{proof}

Note that, by symmetry, a result similar to Lemma \ref{PW.sec.2.lem.5} can be stated for the set $\mathcal{U}$ of the upper tetrahedra of $K$ by reversing the path $(T_i)_{0\leq{i}\leq{k}}$ (or equivalently, by exchanging the below and above relations) while exchanging $\alpha$ with $\gamma$ and $a$ with $c$. In this case, the tetrahedra contained in $\mathcal{P}\cap\mathcal{U}$ must be pulled to $a$ instead of $c$, resulting in a blow-up triangulation $U$ such that $\pi(U^-)$ and $\pi(U^+)$ are equal to $T^\circ$ and $\pi(K^+)$ respectively. It turns out that $L$ and $U$ are not built from the same tetrahedra of $K$.

\begin{lem}\label{PW.sec.2.lem.6}
The sets $\mathcal{L}$ and $\mathcal{U}$ are disjoint.
\end{lem}
\begin{proof}
Consider a tetrahedron $\sigma$ of $K$ contained in $\mathcal{L}$. According to Definition~\ref{PW.sec.2.defn.1}, $K$ must contain a facet $\tau^-$ of $\sigma$ that is disjoint from $\mathrm{P}^\alpha\mathord{\setminus}\{b\}$, not above $\Pi^\beta$, and above $\sigma$ in case it is contained in $\Pi^\beta$. Consider a facet $\tau^+$ of $\sigma$ in $K$ (possibly $\tau^-$ itself). Assume that $\tau^+$ is disjoint from $\mathrm{P}^\gamma\mathord{\setminus}\{b\}$ and not below $\Pi^\beta$. We shall prove that $\tau^+$ is then necessarily contained in $\Pi^\beta$ and above $\sigma$, showing that $\sigma$ cannot be in $\mathcal{U}$ by Definition~\ref{PW.sec.2.defn.1}.

Denote by $\delta$ an arc shared by $\tau^-$ and $\tau^+$. There may be three such arcs if $\tau^+$ coincides with $\tau^-$ and in this case we pick any of them for $\delta$. As $\delta$ is disjoint from $\mathrm{P}^\alpha\mathord{\setminus}\{b\}$ and from $\mathrm{P}^\gamma\mathord{\setminus}\{b\}$, either its two vertices are in $\mathrm{P^\beta}\mathord{\setminus}\pi(\beta)$ or just one of them is in $\mathrm{P^\beta}\mathord{\setminus}\pi(\beta)$ and its other vertex is $b$. The latter situation is impossible because in that case, $\delta$ (and therefore both $\tau^-$ and $\tau^+$) would be above or below $\Pi^\beta$. For the same reason, in the former situation, $\delta$ must be contained in $\Pi^\beta\mathord{\setminus}\beta$. It follows that $\tau^-$ and $\tau^+$ must be above, below, or contained in $\Pi^\beta$. As $\tau^-$ is not above $\Pi^\beta$, it must be below $\Pi^\beta$ or contained in it and in both cases, $\sigma$ is below $\Pi^\beta$ because of the assumption that, if  $\tau^-$ is contained in $\Pi^\beta$, then it is above $\sigma$. Likewise, as $\tau^+$ is not below $\Pi^\beta$, it must be above $\Pi^\beta$ or contained in it. The former case is impossible because $\sigma$  cannot be both below $\Pi^\beta$ and above it. This proves that $\tau^+$ is contained in $\Pi^\beta$ and above $\sigma$, as desired.
\end{proof}

We now prove the following theorem as a consequence of Lemmas \ref{PW.sec.2.lem.5} and~\ref{PW.sec.2.lem.6}. This theorem provides a blow-up triangulation of $\mathrm{P}$ with the same boundary than $K$ and that contains $|\mathcal{L}\cup\mathcal{U}|$ tetrahedra.

\begin{thm}\label{PW.sec.2.thm.1}
There exists a blow-up triangulation $N$ of $\mathrm{P}$ such that
\begin{itemize}
\item[(i)] the number of tetrahedra in $N$ is $|\mathcal{L}\cup\mathcal{U}|$,
\item[(ii)] $\pi(N^-)$ is equal to $\pi(K^-)$ and $\pi(N^+)$ to $\pi(K^+)$, and
\item[(iii)] the tetrahedra of $N$ are obtained from those in $\mathcal{L}\cup\mathcal{U}$ by pulling to $c$ the tetrahedra in $\mathcal{P}\cap\mathcal{L}$ and to $a$ the tetrahedra in $\mathcal{P}\cap\mathcal{U}$.
\end{itemize}
\end{thm}
\begin{proof}
Consider the blow-up triangulation $L$ provided by Lemma \ref{PW.sec.2.lem.5}, whose set of tetrahedra is obtained from $\mathcal{L}$ by pulling to $c$ the tetrahedra that also belong to $\mathcal{P}$. That triangulation is such that $\pi(L^-)$ is equal to $\pi(K^-)$ and $\pi(L^+)$ to $T^\circ$. By symmetry, Lemma \ref{PW.sec.2.lem.5} also provides a blow-up triangulation $U$ whose set of tetrahedra is obtained from $\mathcal{U}$ by pulling to $a$ the tetrahedra contained in $\mathcal{P}\cap\mathcal{U}$. Moreover, $\pi(U^-)$ is equal to $T^\circ$ and $\pi(U^+)$ to $\pi(K^+)$. It follows that $\pi(L^+)$ coincides with $\pi(U^-)$, and one can build the desired triangulation $N$ by gluing $U$ on top of $L$. According to Lemma~\ref{PW.sec.2.lem.6}, the number of tetrahedra in $N$ is $|\mathcal{L}\cup\mathcal{U}|$, as desired.
\end{proof}

\begin{rem}\label{PW.sec.2.rem.1}
Note here that the blow-up triangulation $N$ provided by Theorem \ref{PW.sec.2.thm.1} can be identical to $K$ even when $C$ is not the boundary of a triangle of $K$. For instance, consider the path in the flip-graph of a pentagon depicted in Fig. \ref{PW.sec.1.1.fig.1} and assume that $\alpha$ is the blue dotted arc while $\gamma$ is the red dotted arc. If $K$ is the blow-up triangulation associated to that path, then $\mathcal{P}$ is empty, $\mathcal{L}=\{\psi(3),\psi(4)\}$, and $\mathcal{U}=\{\psi(1),\psi(2)\}$.
\end{rem}

When $X$ is the vertex set of $\mathrm{P}$, we prove the following lemma as a consequence of Theorem \ref{PW.sec.2.thm.1}. It can be thought of as a decomposition lemma because it provides the existence of triangles along which $K$ can be cut, thus making it possible to decompose a blow-up triangulation into smaller blow-up triangulations as highlighted at the beginning of Section \ref{PW.sec.1}.

\begin{lem}\label{PW.sec.2.thm.2}
Consider a convex Euclidean polygon $\mathrm{P}$ with vertex set $X$ and a blow-up triangulation $K$ of $\mathrm{P}$ whose vertex set is also $X$. Further consider a topological circle $C$ obtained as the union of three arcs of $K$. If the number of tetrahedra in $K$ is equal to the distance between $\pi(K^-)$ and $\pi(K^+)$ in $\mathcal{F}(\mathrm{P},X)$, then there exists an arc $\beta$ of $K$ contained in $C$ and a triangle $\tau$ in $S_K(\beta)$ such that $\pi(C)$ is the boundary of $\pi(\tau)$.
\end{lem}
\begin{proof}
Consider a path $(T_i)_{0\leq{i}\leq{k}}$ in $\mathcal{F}(\mathrm{P},X)$ associated to $K$. Assuming that the number $k$ of tetrahedra in $K$ is equal to the distance between $\pi(K^-)$ and $\pi(K^+)$ in $\mathcal{F}(\mathrm{P},X)$, it follows that this path is a geodesic. We assume without loss of generality that one of the arcs of $K$ contained in $C$ belongs to $K^-$. It suffices to discard from the beginning of the path $(T_i)_{0\leq{i}\leq{k}}$ any triangulation that does not contain the image by $\pi$ of such an arc (and to remove from $K$ the corresponding faces). Let us denote that arc by $\alpha$. Likewise, we can assume that an arc $\gamma$ of $K^+$ distinct from $\alpha$ is contained in $C$ by if needed discarding triangulations at the end of the considered path as well as the corresponding faces in $K$.

Denote by $\beta$ the arc of $K$ contained in $C$ and distinct from $\alpha$ and $\gamma$. Assume for contradiction that $\pi(C)$ is not the boundary of the image by $\pi$ of any triangle in $S_K(\beta)$. By Theorem~\ref{PW.sec.2.thm.1}, there exists a blow-up triangulation $N$ of $\mathrm{P}$ with vertex set $X$ whose set of tetrahedra is obtained from $\mathcal{L}\cup\mathcal{U}$ by pulling to $c$ the tetrahedra in $\mathcal{P}\cap\mathcal{L}$ and to $a$ the tetrahedra in $\mathcal{P}\cap\mathcal{U}$. In addition, $\pi(N^-)$ coincides with $\pi(K^-)$ and $\pi(N^+)$ with $\pi(K^+)$. The number of tetrahedra in $N$ is $|\mathcal{L}\cup\mathcal{U}|$ and in particular, $N$ has at most as many tetrahedra as $K$. We shall reach a contradiction by showing that some tetrahedron of $K$ is not contained in $\mathcal{L}\cup\mathcal{U}$. Indeed, in that case, any path from $\pi(K^-)$ to $\pi(K^+)$ in $\mathcal{F}(\mathrm{P},X)$ associated to the blow-up triangulation $N$ is shorter than the distance of $\pi(K^-)$ and $\pi(K^+)$ in $\mathcal{F}(\mathrm{P},X)$.

Since $\pi(C)$ is not the boundary of the image by $\pi$ of a triangle in $S_K(\beta)$, it follows from Lemma~\ref{PW.sec.1.1.lem.2}, that some tetrahedron $\sigma$ in $\mathcal{P}$ is incident to $\beta$. Note that the arc $\delta$ of $K$ that penetrates $C$ and is incident to $\sigma$ must be disjoint from $\beta$. Consider a facet $\tau$ of $\sigma$ in $K$ that is disjoint from $\mathrm{P}^\alpha\mathord{\setminus\{b\}}$ and not above $\Pi^\beta$. As $\tau$ is disjoint from $\mathrm{P}^\alpha\mathord{\setminus\{b\}}$ it cannot be incident to $c$ and therefore, it admits $\delta$ as an edge. As $\tau$ is not above $\Pi^\beta$, the only way for $\varepsilon$ to penetrate $C$ is to be above $\alpha$, which is impossible because $\tau$ is disjoint from $\mathrm{P}^\alpha\mathord{\setminus\{b\}}$. Therefore by Definition~\ref{PW.sec.2.defn.1}, $\sigma$ cannot belong to $\mathcal{L}$. The symmetric argument shows that $\sigma$ cannot belong to $\mathcal{U}$ either.
\end{proof}

\begin{rem}\label{PW.sec.2.rem.2}
The transitive closure of the aboveness relation is a partial order on the faces of a blow-up triangulation. This makes it possible to characterize which of the three arcs contained in $C$ can be taken for $\beta$ in Lemma~\ref{PW.sec.2.thm.2}: one can take for $\beta$ any arc contained in $C$ that is neither below nor above both of the other two with respect to this partial order.
\end{rem}

Let us turn our attention to the case when $X$ contains flat vertices of $\mathrm{P}$. In that case, the proofs in this section all work under the condition that the projections $\langle{T_i}\rangle_\alpha^c$, $\langle\langle{T_i}\rangle_\alpha^c\rangle_\beta^c$, $\langle{T_i}\rangle_\gamma^a$, and $\langle\langle{T_i}\rangle_\gamma^a\rangle_\beta^a$ we used along the way still exist. Indeed, as we discussed in Section \ref{PW.sec.1.5}, the presence of flat vertices may make it impossible to use the projection from \cite{SleatorTarjanThurston1988} and the same holds for the projections we defined in this section: when pulling an arc of a triangulation to a vertex, we must make sure that the pulled arc will not contain a flat vertex in its interior. We can find conditions on the placement of the flat vertices under which the desired projections can be built. As this placement is relative to the arc $\beta$, we need to specify how that arc is chosen.

\begin{lem}\label{PW.sec.2.thm.3}
Consider a convex Euclidean polygon $\mathrm{P}$ and a set $X$ made of all the vertices of $\mathrm{P}$ and finitely-many flat vertices of $\mathrm{P}$. Further consider a blow-up triangulation $K$ of $\mathrm{P}$ with vertex set $X$ and a topological circle $C$ obtained as the union of three arcs $\alpha$, $\beta$, and $\gamma$ of $K$. Assume that
\begin{enumerate}
\item[(i)] $\beta$ is not above both $\alpha$ and $\gamma$ and not below both $\alpha$ and $\gamma$ for the transitive closure of the aboveness relation,
\item[(ii)] if either $\pi(\alpha)$ or $\pi(\gamma)$ is contained in the boundary of $\mathrm{P}$, then their common vertex is not a flat vertex of $\mathrm{P}$, and
\item[(iii)] no flat vertex of $\mathrm{P}$ in $X\cap\mathrm{P}^\beta$ is consecutive to a vertex of $\beta$.
\end{enumerate}

If in addition, the number of tetrahedra contained in $K$ is equal to the distance between $\pi(K^-)$ and $\pi(K^+)$ in $\mathcal{F}(\mathrm{P},X)$, then there exists a triangle $\tau$ in $S_K(\beta)$ such that $\pi(C)$ is the boundary of $\pi(\tau)$.
\end{lem}
\begin{proof}
Denote by $a$, $b$, and $c$ the vertices of $K$ in $C$ with the convention that $a$ is not a vertex of $\alpha$, $b$ is not a vertex of $\beta$, and $c$ is not a vertex of $\gamma$. By Theorem~\ref{PW.sec.2.thm.2}, Remark~\ref{PW.sec.2.rem.2}, and the above discussion, it suffices to show that, for every triangulation $T$ of $\mathrm{P}$ with vertex set $X$, the projections $\langle{T}\rangle_\alpha^c$, $\langle\langle{T}\rangle_\alpha^c\rangle_\beta^c$, $\langle{T}\rangle_\gamma^a$, and $\langle\langle{T}\rangle_\gamma^a\rangle_\beta^a$ exist. By the symmetry that exchanges $\alpha$ with $\gamma$ and $a$ with $c$, we only need to show that $\langle{T}\rangle_\alpha^c$ and $\langle\langle{T}\rangle_\alpha^c\rangle_\beta^c$ exist.
\begin{figure}[b]
\begin{centering}
\includegraphics[scale=1]{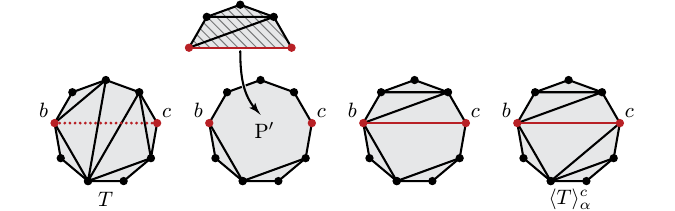}
\caption{An alternative construction for $\langle{T}\rangle_\alpha^c$.}\label{PW.sec.2.fig.3}
\end{centering}
\end{figure}

First recall that if $\pi(\alpha)$ is contained in the boundary of $\mathrm{P}$, then $\mathrm{P}^\alpha$ is equal to $\pi(\alpha)$ and $[T]_\alpha^c$ is equal  to $T$. In that case, $\langle{T}\rangle_\alpha^c$ exists and is equal to $T$ as well by Definition~\ref{PW.sec.2.defn.0}. If, however $\pi(\alpha)$ is not contained in the boundary of $\mathrm{P}$, then consider the polygon $\mathrm{P}'$ obtained as the union of the triangles of $T$ whose interior is non-disjoint from $\mathrm{P}^\alpha$. We shall use the following alternative construction of $\langle{T}\rangle_\alpha^c$ that does not use $[T]_\alpha^c$. First remove from $T$ the arcs contained in $\mathrm{P}'$ and replace them with the image by $\pi$ of the arcs of $K$ contained in $\Pi^\alpha$ as sketched in Fig.~\ref{PW.sec.2.fig.3}. Note that this can always be done, even when $X$ contains flat vertices in $\mathrm{P}^\alpha$ that would prevent $[T]_\alpha^c$ from existing. If $(\mathrm{P}'\mathord{\setminus}\mathrm{P}^\alpha)\cup\pi(\alpha)$ is not just a line segment, we further need to add the arcs of the comb triangulation of that polygon at $c$. The only condition for this comb triangulation to exist is that its arcs should not contain a flat vertex from $X$ in their interior. For that, it suffices to require that $X\cap\mathrm{P}^\beta$ does not contain a flat vertex of $\mathrm{P}$ consecutive to $c$.

Using the same alternative construction to build $\langle\langle{T}\rangle_\alpha^c\rangle_\beta^c$ from $\langle{T}\rangle_\alpha^c$, one can see that this second projection will exist as soon as $\alpha$ is not in the boundary of $\mathrm{P}$ or $b$ is not a flat vertex of $\mathrm{P}$, as desired.
\end{proof}

\section{Several variants of the decomposition lemma}\label{PW.sec.2.5}

The following variant of Lemma~\ref{PW.sec.2.thm.3} states that $C$ is the boundary of the triangle $\tau$ (and not just $\pi(C)$ the boundary of $\pi(\tau)$) under an additional condition on the placement of the flat vertices of $\mathrm{P}$ contained in $X$.

\begin{thm}\label{PW.sec.2.5.thm.1}
Consider a convex Euclidean polygon $\mathrm{P}$, a blow-up triangulation $K$ of $\mathrm{P}$ with vertex set $X$, and three arcs $\alpha$, $\beta$, and $\gamma$ of $K$ whose union is a topological circle $C$. Assume that all of these objects satisfy the same conditions as in the statement of Lemma \ref{PW.sec.2.thm.3}. If in addition, the vertex shared by $\alpha$ and $\gamma$ is not consecutive to a flat vertex of $\mathrm{P}$ contained in $X$, then $C$ is the boundary of a triangle contained in $K$.
\end{thm}
\begin{proof}
According to Lemma~\ref{PW.sec.2.thm.3}, $\pi(C)$ is the boundary of the image by $\pi$ of a triangle from $S_K(\beta)$. Hence, it suffices to show that $K$ does not contain any two arcs with the same pair of endpoints as $\alpha$ or as $\gamma$.

If the vertex shared by $\alpha$ and $\gamma$ is not consecutive to a flat vertex of $\mathrm{P}$ contained in $X$, then, by Lemma \ref{PW.sec.1.5.lem.3}, $\mathcal{F}_\alpha(\mathrm{P},X)$ and $\mathcal{F}_\gamma(\mathrm{P},X)$ are strongly convex subgraphs of $\mathcal{F}(\mathrm{P},X)$. As the number of tetrahedra in $K$ coincides with the distance between $\pi(K^-)$ and $\pi(K^+)$ in $\mathcal{F}(\mathrm{P},X)$, there cannot be two arcs in $K$ with the same pair of endpoints as $\alpha$ or as $\gamma$.
\end{proof}

Let us come back to the important special case when $X$ does not contain any flat vertex of $\mathrm{P}$. In this case, we have proved Lemma \ref{PW.sec.2.thm.3} whereby $\pi(C)$ is the boundary of the image by $\pi$ of a triangle of $K$ incident to $\beta$. Thanks to Theorem \ref{PW.sec.2.5.thm.1}, we can now prove that there is a single such triangle and that $C$ is the boundary of that triangle. As $X$ does not contain any flat vertex, we do not need to care about their placement and can handle at once all the topological circles formed from three arcs of $K$.

We recall that a simplicial complex is \emph{flag} when a subset of its vertices that induces a complete graph in its $1$-skeleton is always the vertex set of one of its faces. It is a consequence of Lemma \ref{PW.sec.1.5.lem.1} that if $X$ is the vertex set of $\mathrm{P}$, then a blow-up triangulation $K$ of $\mathrm{P}$ built from a geodesic path in $\mathcal{F}(\mathrm{P},X)$ is a simplicial complex (which, by Remark \ref{PW.sec.1.rem.2}, is not true in general). According to Theorem \ref{PW.sec.0.thm.0}, this simplicial complex is flag.

\begin{proof}[Proof of Theorem \ref{PW.sec.0.thm.0}]
Assume that $K$ is a blow-up triangulation associated to a geodesic path in $\mathcal{F}(\mathrm{P},X)$. According to Lemma \ref{PW.sec.1.5.lem.1}, $K$ is a simplicial complex. As $X$ does not contain any flat vertex of $\mathrm{P}$, it follows from Theorem~\ref{PW.sec.2.5.thm.1} that any topological circle formed as the union of three arcs of $K$ bounds a triangle of $K$. There remains to show that four vertices of $\mathrm{P}$ such that any two of them are connected by an arc of $K$ form the vertex set of a tetrahedron of $K$. Consider four such vertices. By the above, there exists a topological sphere $\mathcal{S}$ containing all the arcs between two of these vertices as well as four triangles of $K$ each bounded by three of these arcs. In order for the topological ball bounded by $\mathcal{S}$ not to be a tetrahedron of $K$, the interior of that ball must contain at least one vertex of $\mathrm{P}$ which is impossible because the vertices of $\mathrm{P}$ are in the boundary of $K$.
\end{proof}

We are finally ready to prove Theorem \ref{PW.sec.0.thm.1}.

\begin{proof}[Proof of Theorem \ref{PW.sec.0.thm.1}]
Consider a geodesic path in $\mathcal{F}(\mathrm{P},X)$ where $\mathrm{P}$ is a convex polygon and $X$ its vertex set. Further consider the blow-up triangulation $K$ of $\mathrm{P}$ associated to that path. According to Theorem \ref{PW.sec.0.thm.0}, $K$ is flag and as a consequence, if the three edges of a triangle appear in (possibly different) triangulations along the considered path, then that triangle must be contained in some triangulation along that path, as desired.
\end{proof}

\section{Strong convexity fails with two punctures}\label{PW.sec.3}

\begin{figure}[b]
\begin{centering}
\includegraphics[scale=1]{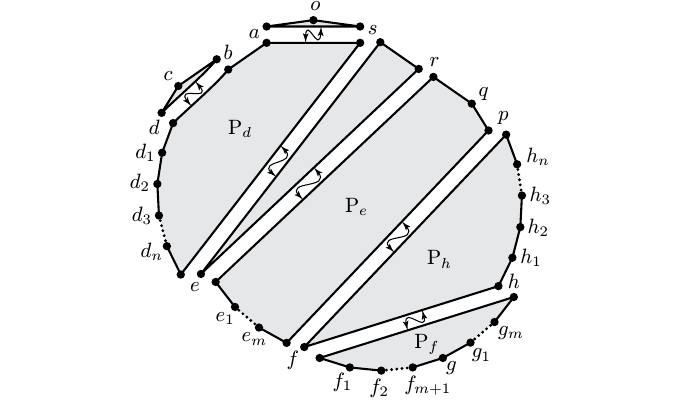}
\caption{The convex polygon $\mathrm{P}$.}\label{PW.sec.3.fig.1}
\end{centering}
\end{figure}

In this section, we consider two positive integers $n$ and $m$, a convex polygon $\mathrm{P}$ with $2n+3m+11$ vertices, and a set $X$ made of all the vertices of $\mathrm{P}$ together with two additional flat vertices. That polygon is sketched and the labels of its vertices are shown on Fig. \ref{PW.sec.3.fig.1}. It is obtained by gluing three triangles together with four convex polygons that we denote by $\mathrm{P}_d$, $\mathrm{P}_e$, $\mathrm{P}_f$, and $\mathrm{P}_h$. Observe that $\mathrm{P}_d$ has a sequence of $n$ vertices labeled $d_1$ to $d_n$ counter-clockwise, $\mathrm{P}_e$ a sequence of $m$ vertices labeled $e_1$ to $e_m$ counter-clockwise, and $\mathrm{P}_h$ a sequence of $n$ vertices labeled $h_1$ to $h_n$ counter-clockwise. Further observe that $\mathrm{P}_f$ has a sequence of $m+1$ vertices labeled $f_1$ to $f_{m+1}$ counter-clockwise, followed by a sequence of $m$ vertices labeled $g_1$ to $g_m$ counter-clockwise. The two flat vertices of $\mathrm{P}$ contained in $X$ are the points denoted by $b$ and $r$ in the figure. In particular, the point $q$, $r$, and $s$ are collinear as well as the points $a$, $b$, and $c$.

\begin{figure}
\begin{centering}
\includegraphics[scale=1]{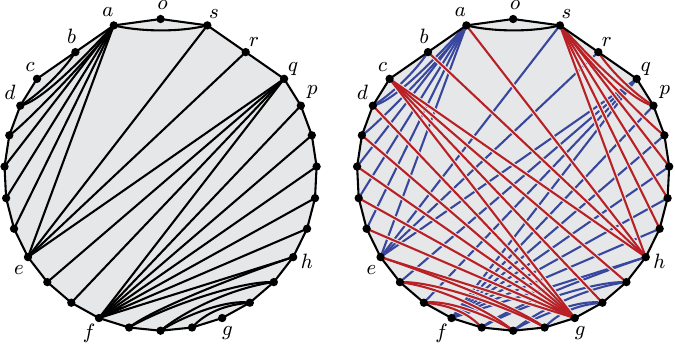}
\caption{The triangulations $T^-$ (left) and $T^+$ (right, on top).}\label{PW.sec.3.fig.2}
\end{centering}
\end{figure}

We now build a triangulation of $\mathrm{P}$ with vertex set $X$ as the union of the comb triangulation of $\mathrm{P}_d$ at vertex $a$, the comb triangulation of $\mathrm{P}_e$ at vertex $q$, the comb triangulation of $\mathrm{P}_h$ at vertex $f$, and a zigzag triangulation of $\mathrm{P}_f$. Recall that a zigzag triangulation is a triangulation such that the arcs that are not contained in the boundary of the polygon form a simple path that alternates between left and right turns. Here, the zigzag triangulation of $\mathrm{P}_f$ we consider is the one where the simple path starts at vertex $f_{m+1}$ and does not contain the vertex $g$. In order to complete this triangulation, one needs to add the three triangles shown Fig. \ref{PW.sec.3.fig.1} as well as their edges. We denote by $T^-$ the resulting triangulation of $\mathrm{P}$. It is depicted on the left of Fig. \ref{PW.sec.3.fig.2} when $n=4$ and $m=2$. In this figure, some arcs have been slightly bent in order to make it easier to see the triangulation near vertices $o$, $c$, and $g$, but all of them can be realized as straight line segments.

We can define a triangulation $T^+$ of $\mathrm{P}$ that is symmetric with $T^-$ with respect to the vertical axis in the representation of Fig. \ref{PW.sec.3.fig.2}. That triangulation is sketched on the right of the figure on top of $T^-$. Observe that the triangulations $T^-$ and $T^+$ share an arc $\eta$ with vertices $a$ and $s$. That arc is precisely the kind of arcs for which Lemma \ref{PW.sec.1.5.lem.3} fails to establish the strong convexity of $\mathcal{F}_\eta(\mathrm{P},X)$ as a subgraph of $\mathcal{F}(\mathrm{P},X)$. We will prove that, in fact, $\mathcal{F}_\eta(\mathrm{P},X)$ is not a strongly convex subgraph of $\mathcal{F}(\mathrm{P},X)$. The strategy is to first establish an upper bound on the distance between $T^-$ and $T^+$ in $\mathcal{F}(\mathrm{P},X)$ and then, to show that their distance in $\mathcal{F}_\eta(\mathrm{P},X)$ is greater than this upper bound when $n$ and $m$ are large enough.

\begin{figure}
\begin{centering}
\includegraphics[scale=1]{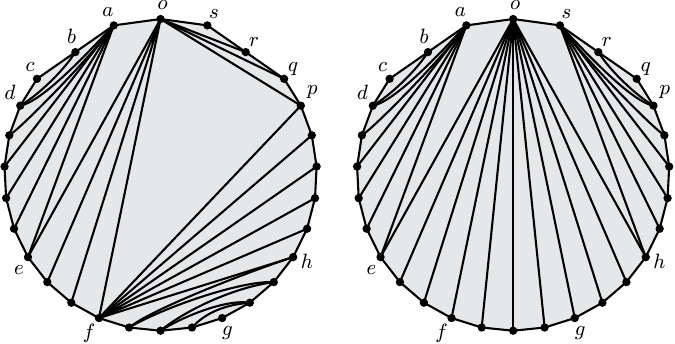}
\caption{The triangulations in the proof of Proposition \ref{PW.sec.3.prop.1}.}\label{PW.sec.3.fig.3}
\end{centering}
\end{figure}

\begin{prop}\label{PW.sec.3.prop.1}
The distance of $T^-$ and $T^+$ in $\mathcal{F}(\mathrm{P},X)$ is at most
$$
2n+6m+24\mbox{.}
$$
\end{prop}
\begin{proof}
We only need to show that the distance in $\mathcal{F}(\mathrm{P},X)$ between $T^-$ and the triangulation $T$ of $\mathrm{P}$ sketched on the right of Fig. \ref{PW.sec.3.fig.3} is at most $n+3m+12$, and the result will follow by symmetry.

We can transform $T^-$ into the triangulation shown on the left of Fig. \ref{PW.sec.3.fig.3} by first performing a sequence of $m+5$ flips that each introduce an arc incident to $o$. From that triangulation, we perform the flip that removes the arc with vertices $o$ and $q$ followed by the flip that removes the arc with vertices $o$ and $r$. Note that both of these flips introduce arcs of $T$. Moreover, the latter introduced arc is one of the arcs of $T$ incident to $s$.

We then perform the flip that replaces the arc with vertices $o$ and $p$ by the arc with vertices $f$ and $s$ followed by a sequence of flips that each introduce one of the remaining $n+1$ arcs of $T$ incident to $s$ which amounts to $n+m+9$ flips so far. From there, $T$ can be reached with a sequence of flips that each introduce one of the remaining $2m+3$ arcs incident to $o$.
\end{proof}

In order to show that the distance between $T^-$ and $T^+$ in $\mathcal{F}_\eta(\mathrm{P},X)$ is greater than $2n+6m+24$ for some choice of $n$ and $m$, we fix from now on a geodesic path $(T_i)_{0\leq{i}\leq{k}}$ from $T^-$ to $T^+$ in $\mathcal{F}_\eta(\mathrm{P},X)$ and the associated blow-up triangulation $K$ of $\mathrm{P}$. Note that the union of the faces of $K$ is a $3$\nobreakdash-dimensional topological ball to which a triangle with vertices $a$, $s$ and $o$ has been glued along an arc whose image by $\pi$ is $\eta$.

We will use a variant of our decomposition lemma that allows us to treat the case when $X$ contains flat vertices of $\mathrm{P}$. This will also illustrate why we call this result a decomposition lemma. In order to be able to use that lemma, we first need to establish the existence of certain arcs of $K$.

\begin{lem}\label{PW.sec.3.lem.1}
If $k$ is less than $3n+6$, then $K$ contains an arc with vertices $a$ and $g$ and, below this arc, an arc with vertices $f$ and $s$.
\end{lem}
\begin{proof}
Denote by $\mathcal{S}^-_x$ the set of the tetrahedra in $K$ incident to and above an arc of $K^-$ that admits $x$ as a vertex. Note that $\mathcal{S}^-_a$ contains $n+2$ different tetrahedra because there are $n+2$ arcs in $K^-$ that are incident to $a$ and not contained in $K^+$ (whose other vertex is $d$, $e$, or $d_i$ where $1\leq{i}\leq{n}$). This follows from the observation that no tetrahedron of $K$ can be incident to two different arcs of $K^-$ and be above them. Likewise, denote by $\mathcal{S}^+_x$ the set of the tetrahedra in $K$ incident to and below an arc of $K^+$ that admits $x$ as a vertex. Observe that again, $\mathcal{S}^+_g$ and $\mathcal{S}^+_s$ each contain $n+2$ different tetrahedra. It turns out that $\mathcal{S}^+_s$ is disjoint from both $\mathcal{S}^-_a$ and $\mathcal{S}^+_g$. Indeed, a tetrahedron in $\mathcal{S}^+_s$ cannot have two vertices on the same side of $\delta$ as $o$, where $\delta$ is the Euclidean line segment between $s$ and $h$. In particular, if $\mathcal{S}^-_a$ and $\mathcal{S}^+_g$ were also disjoint, then the number $k$ of tetrahedra in $K$ would be at least $3n+6$. As a consequence, assuming that $k$ is less than $3n+6$, we can consider a tetrahedron $\sigma^+$ of $K$ such that
$$
\sigma^+\in\mathcal{S}^-_a\cap\mathcal{S}^+_g
$$
and, by symmetry a tetrahedron $\sigma^-$ of $K$ such that
$$
\sigma^-\in\mathcal{S}^+_s\cap\mathcal{S}^-_f\mbox{.}
$$

This immediately implies that $\sigma^-$ has an edge with vertices $f$ and $s$ while $\sigma^+$ has an edge with vertices $a$ and $g$. Finally, observe that $\sigma^-$ is below $\sigma^+$ (otherwise these tetrahedra would not have disjoint interiors) and as a consequence,  the former arc is below the latter, as desired.
\end{proof}
\begin{figure}[b]
\begin{centering}
\includegraphics[scale=1]{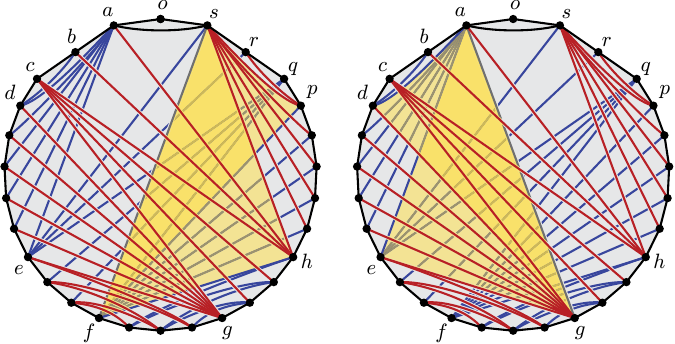}
\caption{Four triangles, two of whose are incident to the arc $\beta^-$ (left) and the two others to the arc $\beta^+$ (right).}\label{PW.sec.3.fig.4}
\end{centering}
\end{figure}

Let us assume that $k$ is less than $3n+6$. In this case, according to Lemma~\ref{PW.sec.3.lem.1}, $K$ must contain two arcs $\beta^-$ and $\beta^+$, the former with vertices $f$ and $s$, and the latter with vertices $a$ and $g$, such that $\beta^-$ is below $\beta^+$. Consider the triangles incident to these arcs shown in Fig.~\ref{PW.sec.3.fig.4}. Two of these triangles on the left of the figure admit $\beta^-$ as an edge, while $\beta^+$ is an edge of the other two on the right of the figure. The other two edges of each of these triangles are an arc contained in $K^-$ and an arc contained in $K^+$. We will denote by $C_d$, $C_e$, $C_h$, and $C_p$ the topological circles bounding these triangles, in such a way that $C_x$ contains $x$. It turns out that we can apply Theorem \ref{PW.sec.2.5.thm.1} for each of these circles within $K$. For instance, the circle $C_p$ is obtained by gluing the arc $\alpha$ of $K^-$ between vertices $f$ and $p$, the arc $\gamma$ of $K^+$ with vertices $s$ and $p$, and the arc $\beta^-$. As $\alpha$ is in $K^-$ and $\gamma$ in $K^+$, it is immediate that $\beta^-$ is not below $\alpha$ and not above $\gamma$ with respect to the transitive closure of the aboveness relation. Moreover, neither $\pi(\alpha)$ nor $\pi(\gamma)$ are contained in the boundary of $\mathrm{P}$. Now recall that the only two flat vertices of $\mathrm{P}$ contained in $X$ are $b$ and $r$. As $b$ is not consecutive to $f$ or $s$ in $X$ while $r$ is not consecutive to $p$, it follows from Theorem~\ref{PW.sec.2.5.thm.1}, that $K$ contains a triangle $\Delta_p$ bounded by $\alpha$, $\beta^-$, and $\gamma$. We identify $\Delta_p$ with the triangle bounded by these arcs shown on the left of Fig.~\ref{PW.sec.3.fig.4}. We can use the same argument for the three other circles $C_d$, $C_e$, and $C_h$: according to Theorem~\ref{PW.sec.2.5.thm.1}, these circles bound three triangles $\Delta_d$, $\Delta_e$, and $\Delta_h$ of $K$, respectively, that we identify with the triangles shown in Fig.~\ref{PW.sec.3.fig.4}.

Now consider the family $(\Pi_i)_{0\leq{i}\leq{k}}$ of topological disks corresponding to $K$ and to $(T_i)_{0\leq{i}\leq{k}}$ via Proposition \ref{PW.sec.1.lem.1}.
\begin{figure}[b]
\begin{centering}
\includegraphics[scale=1]{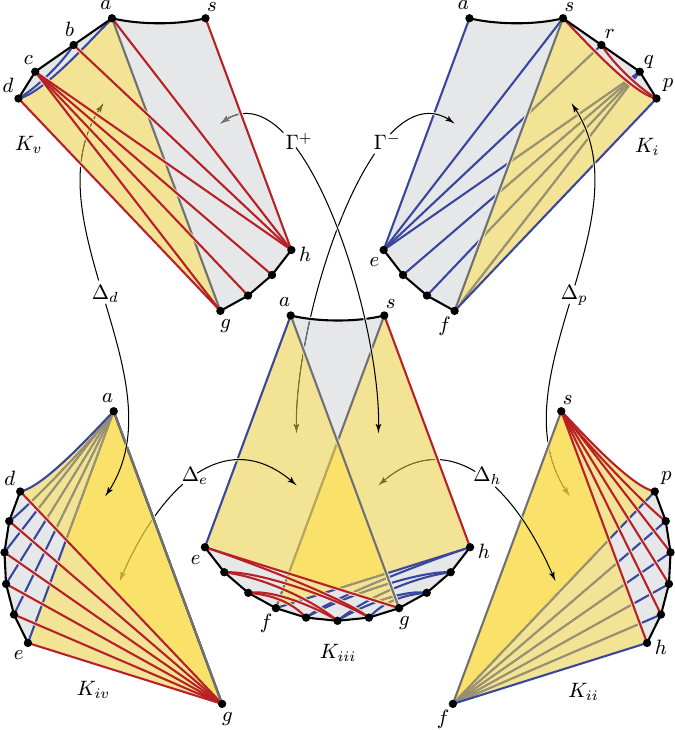}
\caption{A decomposition of $K$.}\label{PW.sec.3.fig.5}
\end{centering}
\end{figure}
Pick two integers $i$ and $j$ such that the triangles $\Delta_h$ and $\Delta_e$ are contained in $\Pi_i$ and $\Pi_j$, respectively. Note that $i$ must be less than $j$ because $\beta^-$ is below $\beta^+$. As $\Delta_e$ is contained in $\Pi_j$ and $i$ is less than $j$, it follows from Proposition \ref{PW.sec.1.lem.1} that the arc $\delta$ in $K^-$ with vertices $a$ and $e$ is necessarily contained in $\Pi_i$. Similarly, the arc $\varepsilon$ in $K^+$ with vertices $s$ and $h$ is contained in $\Pi_j$ because $\Delta_h$ is contained in $\Pi_i$. Denote by $\Gamma^-$ the portion of $\Pi_i$ made of the points whose image by $\pi$ is not on the same side of $\pi(\beta^-)$ as $h$ and not on the same side of $\pi(\delta)$ as $d$. Likewise, let $\Gamma^+$ be the set of all points in $\Pi_j$ whose image by $\pi$ is not on the same side of $\pi(\beta^+)$ as $e$ and not on the same side of $\pi(\varepsilon)$ as $p$.

If follows from Proposition \ref{PW.sec.1.lem.1}, that $\Gamma^-$ and $\Gamma^+$ are unions of triangles of $K$. As in addition the four triangles $\Delta_d$, $\Delta_e$, $\Delta_h$, and $\Delta_p$ belong to $K$, we can cut $K$ along $\Gamma^-$, $\Gamma^+$, and these four triangles. This results in the five blow-up triangulations $K_i$, $K_{ii}$, $K_{iii}$, $K_{iv}$, and $K_v$ of smaller polygons depicted in Fig. \ref{PW.sec.3.fig.5}. It is important to keep in mind that the decomposition of $K$ into $K_i$, $K_{ii}$, $K_{iii}$, $K_{iv}$, and $K_v$ exists thanks to the assumption that $k$ is less $3n+6$ so that we know from Lemma \ref{PW.sec.3.lem.1} that the arcs $\beta^-$ and $\beta^+$ belong to $K$. %Let us assume for a moment that this is the case in order to study $K_{iii}$.
We can estimate the number of tetrahedra contained in $K_{iii}$ using Theorem \ref{PW.sec.1.1.thm.1}. In fact, we will estimate the number of tetrahedra in the blow-up triangulation $\tilde{K}_{iii}$ sketched at the top of Fig. \ref{PW.sec.3.fig.6}, on the left. It is obtained from $K_{iii}$ by removing all the faces incident to both $a$ and $e$ or to both $h$ and $s$, by pulling to $e$ the remaining faces incident to $a$, and by pulling to $h$ the remaining faces incident to $s$. In that process, the faces incident to both $a$ and $s$ (but not to $e$ or $h$) will first be pulled to $e$ and then to $h$. Note that $\tilde{K}_{iii}$ exists by Lemma \ref{PW.sec.1.1.lem.4} and that it contains at most as many tetrahedra as $K_{iii}$. Let us denote by $l^-$ the number of arcs in $K_{iii}^-\mathord{\setminus}K_{iii}^+$ that are incident to $s$ and by $l^+$ the number of arcs in $K_{iii}^+\mathord{\setminus}K_{iii}^-$  that are incident to $a$. By construction, $\tilde{K}_{iii}^-\mathord{\setminus}\tilde{K}_{iii}^+$ contains $l^-+1$ arcs incident to $h$ and $\tilde{K}_{iii}^+\mathord{\setminus}\tilde{K}_{iii}^-$ contains $l^++1$ arcs incident to $e$.

\begin{lem}\label{PW.sec.3.lem.2}
There are at least $5m+4-l^--l^+$ tetrahedra in $\tilde{K}_{iii}$.
\end{lem}
\begin{proof}
Consider the triangle in $\tilde{K}_{iii}^+$ incident to $e$ and $h$ and denote by $\delta$ the edge of that triangle that is not incident to $e$. If the vertices of $\delta$ are $g_m$ and $h$, then denote by $\Omega$ the convex polygon whose vertices are the vertices of $\tilde{K}_{iii}$ and by $\Xi$ the Euclidean line segment with vertices $g_m$ and $h$. Otherwise, consider the convex hull of the vertices of $\tilde{K}_{iii}$, cut that convex hull along $\pi(\delta)$, and denote by $\Omega$ and $\Xi$ the two polygons resulting from that cut with the convention that $\Omega$ contains $e$ and $\Xi$ contains $g_m$. Further denote by $x$ the vertex of $\delta$ other than $h$. By construction, the conditions of Theorem~\ref{PW.sec.1.1.thm.1} are satisfied by $\tilde{K}_{iii}$ with respect to $\delta$, $\Xi$ and $\Omega$. Note in particular that the condition regarding flat vertices is void because no vertex of $\tilde{K}_{iii}$ is a flat vertex of $\Omega\cup\Xi$. Moreover, two edges of $\Omega\cup\Xi$ contained in $\Xi$ have distinct links in $\tilde{K}_{iii}^-$ and these links are contained in $\Omega$ as shown in Fig.~\ref{PW.sec.3.fig.6}. These links are also disjoint from the link of $\delta$ in $\tilde{K}_{iii}^+$ (because the point of $\Omega$ contained in the latter link is $e$). Now consider the blow-up triangulation $L$ provided by Theorem \ref{PW.sec.1.1.thm.1} for this choice of $\delta$, $\Xi$ and $\Omega$: it is obtained from $\tilde{K}_{iii}$ by removing the faces incident to more than one vertex of $\Xi$ and by pulling to $x$ the remaining faces incident to a point in $\Xi\mathord{\setminus}\{x\}$ as sketched at the top of Fig. \ref{PW.sec.3.fig.6}, on the right. 
\begin{figure}
\begin{centering}
\includegraphics[scale=1]{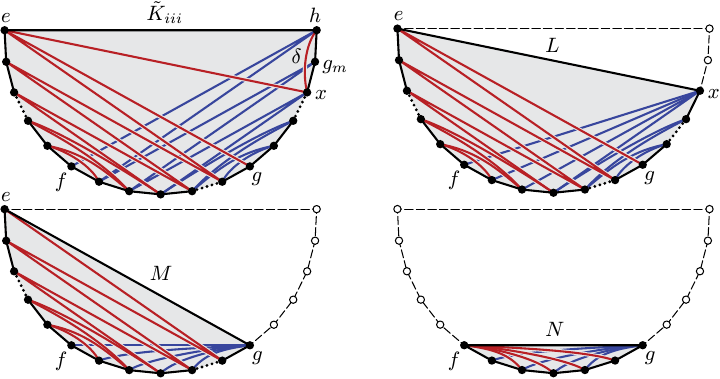}
\caption{The four blow-up triangulations $\tilde{K}_{iii}$, $L$, $M$, and $N$ used in the proof of Lemma \ref{PW.sec.3.lem.2}.}\label{PW.sec.3.fig.6}
\end{centering}
\end{figure}
Moreover the number of tetrahedra in $L$ is less than the number of tetrahedra in $K_{iii}$ by at least $2|X\cap\Xi|-3$. It is important to note here that the number of arcs incident to $e$ in $L^+\mathord{\setminus}L^-$ is now $l^+$ (so this number is down by one compared to $\tilde{K}_{iii}$).

Now observe that, if $x$ is not equal to $g$, then one can invoke Theorem~\ref{PW.sec.1.1.thm.1} again, this time by considering the triangle of $L^+$ incident to $e$ and $x$, by using for $\delta$ the edge of that triangle that does not admit $e$ as a vertex, and for $\Omega$ and $\Xi$ the polygons (or, possibly a line segment for $\Xi$) obtained by cutting the convex hull of $L$ along that arc with the convention that $\Omega$ contains $e$. Repeating this procedure, one eventually ends up with the blow-up triangulation $M$ shown at the bottom of Fig. \ref{PW.sec.3.fig.6}, on the left.

Let us estimate the number of tetrahedra in $M$. By Theorem \ref{PW.sec.1.1.thm.1}, this number at least is twice the number of the vertices of $\tilde{K}_{iii}$ that are no longer vertices of $M$ minus the number of times we invoked that theorem. However, each time we invoked Theorem \ref{PW.sec.1.1.thm.1}, we have lost exactly one of the $l^++1$ arcs incident to $e$ in $\tilde{K}_{iii}^+\mathord{\setminus}\tilde{K}_{iii}^-$ and there is only one such arc left in $M$. As $M$ has $m+1$ vertices less than $\tilde{K}_{iii}$, the number of tetrahedra in $M$ is less than the number of tetrahedra in $\tilde{K}_{iii}$ by at least $2m+2-l^+$.

Now observe that the same procedure can be performed within $M$, but this time by considering the triangle in $M^-$ incident to $e$ and $g$ and its edge $\delta$ that is not incident to $g$. In this case, Theorem \ref{PW.sec.1.1.thm.1} can still be invoked with respect to the arc $\delta$ but one should take care to exchange $M^-$ and $M^+$ in the statement of the theorem. The conditions of this theorem remain satisfied during the process, as we invoke Theorem \ref{PW.sec.1.1.thm.1} for possibly a succession of blow-up triangulations. This process ends up with the triangulation $N$ sketched at the bottom of Fig. \ref{PW.sec.3.fig.6}, on the right. Similarly as between $\tilde{K}_{iii}$ and $M$, the number of tetrahedra in $N$ is less than the number of tetrahedra in $M$ by at least $2m+2-l^-$. The number of tetrahedra in $N$ is at least $m$ as there are $m$ arcs in $N^-\mathord{\setminus}N^+$, which completes the proof.
\end{proof}

We now estimate the number of tetrahedra in $K_v$. By symmetry, this also makes it possible to estimate the number of tetrahedra in $K_i$.

\begin{lem}\label{PW.sec.3.lem.3}
There are at least $m+l^+$ tetrahedra in $K_v$.
\end{lem}
\begin{proof}
Consider the set $\mathcal{S}^+$ of the tetrahedra in $K_v$ incident to and below an arc of $K_v^+$. Looking at the sketch of $K_v$ in Fig. \ref{PW.sec.3.fig.5}, one can see that all the arcs of $K_v^+$ are above an arc of $K_v^-$ except possibly the one incident to $a$ and $h$. As $K_v^+\mathord{\setminus}K_v^-$ contains $m+4$ arcs and a tetrahedron in $K_v$ cannot be incident and below two distinct arcs of $K_v^+$,
\begin{equation}\label{PW.sec.3.lem.3.eq.1}
|\mathcal{S}^+|\geq{m+3}\mbox{.}
\end{equation}

Now let $\mathcal{S}_a^-$ be the set of the tetrahedra from $K_v$ that are incident to and above an arc incident to $a$ and contained in $K_v^-$. As above, at most one of these arcs can be contained in $K_v^+$ and therefore,
\begin{equation}\label{PW.sec.3.lem.3.eq.2}
|\mathcal{S}_a^-|\geq{l^+-1}\mbox{.}
\end{equation}

Only two tetrahedra can belong to both $\mathcal{S}^+$ and $\mathcal{S}_a^-$. Indeed, as one can see in Fig. \ref{PW.sec.3.fig.5} all the tetrahedra in $\mathcal{S}^+$ are incident to $c$ except possibly two. Hence, if there were three or more tetrahedra in $\mathcal{S}^+\cap\mathcal{S}_a^-$, then one of them would be incident to both $a$ and $c$ and its edge with these vertices would contain $b$ in its interior. Hence, $K_v$ contains at least $|\mathcal{S}^+|+|\mathcal{S}_a^-|-2$ tetrahedra. The desired result therefore follows from (\ref{PW.sec.3.lem.3.eq.1}) and (\ref{PW.sec.3.lem.3.eq.2}).
\end{proof}

Combining Lemmas \ref{PW.sec.3.lem.2} and \ref{PW.sec.3.lem.3} provides the following.

\begin{lem}\label{PW.sec.3.lem.4}
The distance of $T^-$ and $T^+$ in $\mathcal{F}_\eta(\mathrm{P},X)$ is at least
$$
\min\{3n+6,2n+7m+6\}\mbox{.}
$$
\end{lem}
\begin{proof}
Assume that the distance of $T^-$ and $T^+$ in $\mathcal{F}_\eta(\mathrm{P},X)$ is less than $3n+6$. One can then decompose $K$ into the blow-up triangulations $K_i$, $K_{ii}$, $K_{iii}$, $K_{iv}$, and $K_v$. Recall that the distance of $T^-$ and $T^+$ in $\mathcal{F}_\eta(\mathrm{P},X)$ is the total number of tetrahedra in these five blow-up triangulations. We shall prove that this number of tetrahedra is at least $2n+7m+6$.

Observe that the number of tetrahedra in $K_{ii}$ is at least $n+1$ (because there are $n+1$ arcs in $K_{ii}^-\mathord{\setminus}K_{ii}^+$). Similarly, the number of tetrahedra in $K_{iv}$ is also at least $n+1$. By Lemma \ref{PW.sec.3.lem.3}, $K_v$ contains at least $m+l^+$ tetrahedra and by symmetry, $K_i$ contains at least $m+l^-$ tetrahedra. Finally, according to Lemma \ref{PW.sec.3.lem.2}, there are at least $5m+4-l^--l^+$ tetrahedra in $\tilde{K}_{iii}$, and the result is obtained by summing these numbers.
\end{proof}

If $n$ is greater than $6m+18$ and $m$ is greater than $18$, then by Proposition~\ref{PW.sec.3.prop.1} and Lemma \ref{PW.sec.3.lem.4}, the distance of $T^-$ and $T^+$ in $\mathcal{F}(\mathrm{P},X)$ is less than their distance in $\mathcal{F}_\eta(\mathrm{P},X)$. In particular, $\mathcal{F}_\eta(\mathrm{P},X)$ is then not a weakly convex subgraph of $\mathcal{F}(\mathrm{P},X)$ because no geodesic path between $T^-$ and $T^+$ is contained in $\mathcal{F}_\eta(\mathrm{P},X)$. It is noteworthy that, by choosing $m$ and $n$ large enough, we can make the gap between the distance of $T^-$ and $T^+$ in $\mathcal{F}_\eta(\mathrm{P},X)$ and their distance in $\mathcal{F}(\mathrm{P},X)$ arbitrarily large. This proves Theorem \ref{PW.sec.0.thm.2} in the case when $X$ contains exactly two flat vertices of $\mathrm{P}$ and no puncture. Let us now prove this theorem in the general case when the number of punctures and flat vertices of $\mathrm{P}$ contained in $X$ sum to at least $2$.

\begin{proof}[Proof of Theorem \ref{PW.sec.0.thm.2}]
As discussed above, the theorem holds when there are two flat vertices of $\mathrm{P}$ in $X$ and no puncture. Let us first extend this to at least three flat vertices and no puncture. Glue a triangle $\tau$ along the edge $\zeta$ of $\mathrm{P}$ with vertices $a$ and $o$ in order to build a polygon $\mathrm{Q}$. We place the vertex $p$ of $\tau$ opposite $\zeta$ close enough to $o$ for $\mathrm{Q}$ to remain convex. We also require $s$, $o$, and $p$ to be collinear in order for $o$ to become a flat vertex of $\mathrm{Q}$. Denote by $Y$ a set made of $p$ and finitely-many flat vertices of $\mathrm{Q}$ contained in the interior of the line segment between $p$ and $o$.

Consider the natural graph isomorphism
$$
\varphi:\mathcal{F}(\mathrm{P},X)\rightarrow\mathcal{F}_\zeta(\mathrm{Q},X\cup{Y})
$$
such that $\varphi(T)$ is the only triangulation of $\mathrm{Q}$ that contains $T$ as a subset. Observe that $\varphi$ sends $\mathcal{F}_\eta(\mathrm{P},X)$ to $\mathcal{F}_\eta(\mathrm{Q},X\cup{Y})$. As $\mathcal{F}_\eta(\mathrm{P},X)$ is not weakly convex in $\mathcal{F}(\mathrm{P},X)$ it follows that $\mathcal{F}_\eta(\mathrm{Q},X\cup{Y})$ is not weakly convex in $\mathcal{F}_\zeta(\mathrm{Q},X\cup{Y})$ and it cannot be weakly convex in $\mathcal{F}(\mathrm{Q},X\cup{Y})$ either.

Now that the theorem has been established in the special case when $X$ does not contain any puncture, the general case follows from the observation that a convex quadrilateral whose vertices are displaced by a sufficiently small amount will remain convex. As a consequence, one can turn into punctures a subset of the flat vertices of $\mathrm{P}$ contained in $X$ without modifying $\mathcal{F}(\mathrm{P},X)$. It suffices to move them into the interior of $\mathrm{P}$ in such a way that they remain close enough to their initial position. As this operation does not affect $\mathcal{F}(\mathrm{P},X)$, this proves the theorem in the general case. 
\end{proof}

\section{Large gaps in distance estimates from arc crossings}\label{PW.sec.4}

It is well-known that, when $\mathrm{P}$ is a convex Euclidean polygon and $X$ denotes its vertex set, the distance in $\mathcal{F}(\mathrm{P},X)$ between any two triangulations is at most $2|X|-10$ when $|X|>12$ \cite{SleatorTarjanThurston1988}. It is also known that this bound is sharp: one can always pick two triangulations whose distance in $\mathcal{F}(\mathrm{P},X)$ is $2|X|-10$ \cite{Pournin2014,SleatorTarjanThurston1988}. However, there is no efficient algorithm for determining the distance in $\mathcal{F}(\mathrm{P},X)$ between any two triangulations \cite{ClearyMaio2018,Dorfer2026,KanjSedgwickXia2017}. Interestingly, computing these distances is an important problem for a number of applications (see \cite{ClearyMaio2018} and references therein). In order to estimate these distances anyways, a popular method is based on the number of arc crossings~\cite{ClearyMaio2018}. Indeed, consider two triangulations $T^-$ and $T^+$ of a topological surface $\Sigma$ with the same vertex set $X$. If $T^-$ does not coincide with $T^+$, one can always find an arc $\varepsilon$ in $T^-$ such that the arc introduced when $\varepsilon$ is flipped in $T^-$ crosses fewer arcs from $T^+$ than $\varepsilon$ does. This has been proven in~\cite{HankeOttmannSchuierer1996} for triangulations of punctured Euclidean polygons and in~\cite{DisarloParlier2019} for triangulations of topological surfaces. Therefore, by performing a flip in $T^-$ such that the number of arc crossings with $T^+$ decreases most, and by repeating this procedure, one eventually reaches $T^+$. The length of the path in $\mathcal{F}(\Sigma,X)$ thus constructed from $T^-$ to $T^+$ can be used as an estimate for the distance between these triangulations in $\mathcal{F}(\Sigma,X)$.

Note that several paths of different lengths may be obtained from this procedure in case different flips in a triangulation make the number of arc crossings with $T^+$ decrease most. We shall denote by $\ell(T^-,T^+)$ the length of the shortest possible path resulting from this method.

Little is known on how accurate the estimate provided by this arc crossings based method. Only recently has it been observed the path in $\mathcal{F}(\Sigma,X)$ obtained from this method is not always a geodesic \cite{ClearyMaio2018}. In other words, this distance estimate is sometimes one off from the distance in $\mathcal{F}(\Sigma,X)$ of the two considered triangulations. In this section, we prove that this method can overestimate the distances in $\mathcal{F}(\Sigma, X)$ by a multiplicative factor arbitrarily close to $3/2$. We first prove this in the special case when $\Sigma$ is a convex polygon, and then when $\Sigma$ is an oriented topological surface.

\begin{figure}[b]
\begin{centering}
\includegraphics[scale=1]{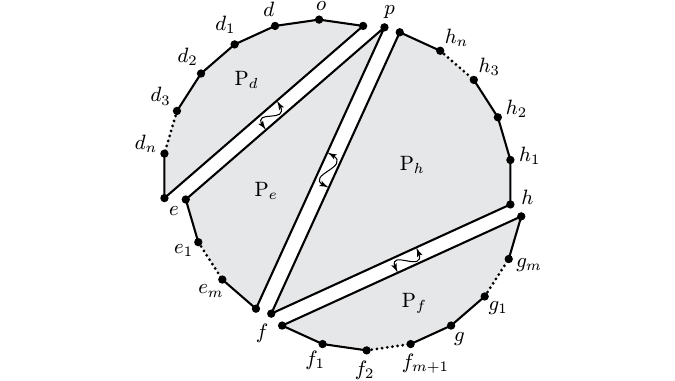}
\caption{The convex polygon $\mathrm{P}$.}\label{PW.sec.4.fig.1}
\end{centering}
\end{figure}

\begin{figure}
\begin{centering}
\includegraphics[scale=1]{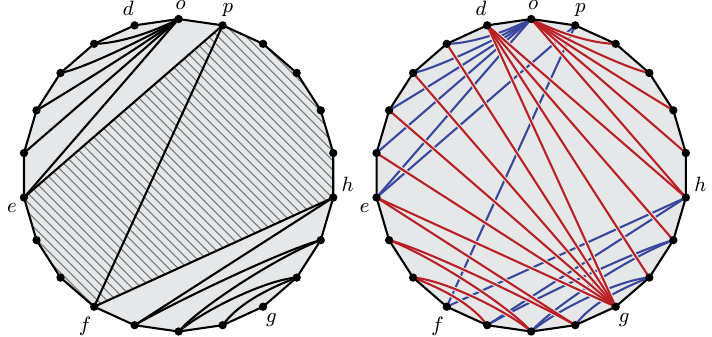}
\caption{The triangulations $T^-$ (left) and $T^+$ (right, on top).}\label{PW.sec.4.fig.2}
\end{centering}
\end{figure}

Let us first describe a family of triangulation pairs of the convex Euclidean polygon $\mathrm{P}$ with $2n+3m+8$ vertices obtained by gluing the four polygons $\mathrm{P}_d$, $\mathrm{P}_e$, $\mathrm{P}_f$, and $\mathrm{P}_h$ as shown in Fig. \ref{PW.sec.4.fig.1}. Note that $\mathrm{P}_f$ and $\mathrm{P}_h$ coincide with the polygons already shown in Fig. \ref{PW.sec.3.fig.1}. 
Further observe that $\mathrm{P}_d$ has a sequence of $n$ vertices labeled $d_1$ to $d_n$ and $\mathrm{P}_e$ a sequence of $m$ vertices labeled $e_1$ to $e_m$ counter-clockwise. Here, $X$ is the vertex set of $\mathrm{P}$ and does not contain any flat vertex or puncture. Now consider a triangulation $T^-$ of $\mathrm{P}$ that contains the arc with vertices $f$ and $p$, the comb triangulation of $\mathrm{P}_d$ at vertex $o$ and the same zigzag triangulation of $\mathrm{P}_f$ as in Section \ref{PW.sec.3}. This triangulation is sketched on the left of Fig. \ref{PW.sec.4.fig.2}. Note that this description does not tell all the arcs of $T^-$. In particular $T^-$ can admit as a subset any triangulation of $\mathrm{P}_e$ and any triangulation of $\mathrm{P}_h$ (these two polygons are striped on the left of Fig. \ref{PW.sec.4.fig.2}). The special case when $T^-$ admits as subsets the comb triangulation of $\mathrm{P}_e$ at vertex $p$ and the comb triangulation of $\mathrm{P}_h$ at vertex $f$ will be particularly interesting. In fact, we will denote by $T^+$ the triangulation of $\mathrm{P}$ that is symmetric to this special case with respect to the vertical axis in the representation of Fig.~\ref{PW.sec.4.fig.2}. The two triangulations are shown together on the right of Fig.~\ref{PW.sec.4.fig.2}, where $T^+$ is shown on top of $T^-$. Observe that, just as with Fig. \ref{PW.sec.3.fig.2} some of the arcs shown in Fig.~\ref{PW.sec.4.fig.2} are bent in order to better show the triangulations near vertices $d$, $f$, $g$ and $p$, but all these arcs can be realized as straight line segments.

Let us now look at how the number of arc crossings with $T^+$ varies when an arc of $T^-$ is flipped. Our goal here is to identify which arc may be flipped first by the method that we described above. Denote by $\mathcal{A}_f$ the set of the arcs in $T^-$ whose two vertices are $f$ and $h_i$ where $1\leq{i}\leq{n}$. Similarly, denote by $\mathcal{A}_p$ the set of the arcs in $T^-$ that are incident to $p$ on one end and to $e_i$ on the other, where $1\leq{i}\leq{m}$. We first show that, when $\mathcal{A}_f\cup\mathcal{A}_p$ contains sufficiently many arcs, flipping any other arc in $T^-$ does not decrease the number of arc crossings with $T^+$ by too much. This is proven by considering several families of arcs from $T^-$.

\begin{prop}\label{PW.sec.4.prop.1}
If $n$ is at least $m$ and $\mathcal{A}_f\cup\mathcal{A}_p$ contains at least $m+3$ arcs, then flipping any of the arcs of $T^-$ that do not belong to $\mathcal{A}_f\cup\mathcal{A}_p$ makes the number of arc crossings with $T^+$ decrease by at most $n+m$. %In addition, there exists an arc in $\pi(\mathcal{A}_f^-\cup\mathcal{A}_p^-)$ whose flip in $T^-$ makes that number decrease by more than $n+m$.
\end{prop}
\begin{proof}
Assume that $\mathcal{A}_f\cup\mathcal{A}_p$ contains at least $m+3$ arcs. By construction, $|\mathcal{A}_f\cup\mathcal{A}_p|$ is at most $n+m$, and therefore, $n$ is at least $3$.

Note that the arc of $T^-$ with vertices $o$ and $d_i$, where $1\leq{i}\leq{n}$ crosses exactly $m+i+1$ arcs of $T^+$ as can be seen in Fig. \ref{PW.sec.4.fig.2}. If that arc is flipped in $T^-$, the arc that replaces it crosses a single arc from $T^+$ and therefore, this flip makes the number of arc crossings with $T^+$ decrease by at most $n+m$. The arc of $T^-$ with vertices $e$ and $o$ crosses $n+m+2$ arcs of $T^+$ and flipping it introduces an arc that crosses $2n+m+2$ arcs of $T^+$. Hence, the desired result holds for this arc as well. Now consider an arc of $T^-\mathord{\setminus}T^+$ incident to $f_i$ where $1\leq{i}\leq{m}$ and observe that, when that arc is flipped in $T^-$, the number of crossings with $T^+$ varies by at most $1$. As $n\geq3$, the desired result holds for any such arc. In addition, flipping the arc of $T^-$ with vertices $f_{m+1}$ and $g_1$ decreases the number of arc crossings with $T^+$ by exactly $n$ and the proposition also holds for that arc.

Let us now check the arcs of $T^-$ contained in $\mathrm{P}_e$ or $\mathrm{P}_h$. Observe that the arc with vertices $e$ and $p$ crosses $2n+m+3$ arcs of $T^+$. Flipping that arc introduces an arc incident to $o$ and whose other vertex is either $f$ or $e_i$ where $1\leq{i}\leq{m}$. The introduced arc crosses at least $n+m+4$ arcs of $T^+$. Thus that flip makes the number of arc crossings decrease by at most $n-1$, as desired. 
%Recall that $|\mathcal{A}_f^-\cup\mathcal{A}_p^-|\geq{m+3}$. By construction, $\mathcal{A}_p^-$ contains at most $m$ arc and as a consequence $\mathcal{A}_f^-$ is non-empty.
The arc of $T^-$ with vertices $f$ and $h$ crosses $n+3m+3$ arcs of $T^+$. When this arc is flipped, it is replaced by an arc $\delta$ with vertices $f_1$ and $h_i$ where $1\leq{i}\leq{n}$ (note in particular that $\delta$ cannot be incident to $p$ because as $\mathcal{A}_f\cup\mathcal{A}_p$ contains more than $m$ arcs, $\mathcal{A}_f$ is not empty and the arcs it contains would then cross $\delta$). The arc $\delta$ crosses at least $n+3m+4$ arcs of $T^+$ and the result therefore also holds in this case.

Now consider the arc $\varepsilon$ of $T^-$ with vertices $f$ and $p$. For convenience, let us denote $e$ by $e_0$. Observe that $\varepsilon$ crosses $2n+3m+5$ arcs of $T^+$ (that is, all the arcs of $T^+$). When flipped, $\varepsilon$ is replaced by an arc $\delta$ incident to $e_i$ with $0\leq{i}\leq{m}$ and $h_j$ with $1\leq{j}\leq{n}$. Note again, that $\delta$ cannot be incident to $h$ because $\mathcal{A}_f$ is non-empty. In particular, $\delta$ crosses exactly $n+m+2i+j+2$ arcs of $T^+$. As $\mathcal{A}_f\cup\mathcal{A}_p$ contains at least $m+3$ arcs, 
$$
i+j\geq{m+3}\mbox{.}
$$

Therefore, $\delta$ crosses at least $n+2m+5$ arcs of $T^+$ and the number of arc crossings decreases by at most $n+m$ when $\varepsilon$ is flipped in $T^-$.

Finally, observe that any of the arcs of $T^-$ that is contained in $\mathrm{P}_e$ but not incident to $p$ crosses at most $2m$ arcs of $T^+$. Similarly, any arc of $T^-$ contained in $\mathrm{P}_h$ but not incident to $f$ crosses at most $n$ arcs of $T^+$. If $m\leq{n}$, then the number of crossing with $T^+$ cannot decrease by more than $n+m$ when such an arc is flipped in $T^-$, which completes the proof.
\end{proof}

Now let us look at the arcs in $\mathcal{A}_f\cup\mathcal{A}_p$.% in $T^-$ decreases the number of arc crossings with $T^+$ by more than $n+m$. %We are particularly interested in the arcs that belong to $\mathcal{A}_f^-$ or to $\mathcal{A}_p^-$

\begin{prop}\label{PW.sec.4.prop.2}
If $n$ is positive, at least $2m$ and at most $m|\mathcal{A}_f\cup\mathcal{A}_p|$, then there exists an arc in $\mathcal{A}_f\cup\mathcal{A}_p$ such that flipping that arc in $T^-$ makes the number of arc crossings with $T^+$ decrease by more than $n+m$.
\end{prop}
\begin{proof}
Assume that $n$ is positive and that
\begin{equation}\label{PW.sec.4.prop.2.eq.1}
2m\leq{n}\leq{m|\mathcal{A}_f\cup\mathcal{A}_p|}\mbox{.}
\end{equation}

First consider the case when $\mathcal{A}_p$ is non-empty and pick an arc in that set. Such an arc crosses at least $2n+m+5$ arcs of $T^+$. When the considered arc is flipped in $T^-$, the arc that replaces it crosses at most $2m$ arcs of $T^+$. Hence, the number of arc crossings with $T^+$ decreases by at least $2n-m+5$ which, by (\ref{PW.sec.4.prop.2.eq.1}) is greater than $n+m$ and the result holds in this case.

Let us now review the case when $\mathcal{A}_p$ is empty. As $n$ is positive, we obtain from (\ref{PW.sec.4.prop.2.eq.1}) that $\mathcal{A}_f$ is non-empty and that $m$ is positive. For each arc $\varepsilon$ contained in $\mathcal{A}_f$, denote by  $u(\varepsilon)$ the number of arcs of $T^+$ crossed by the arc $\delta$ introduced when $\varepsilon$ is flipped within $T^-$. By construction,
\begin{equation}\label{PW.sec.4.prop.2.eq.2}
\sum_{\varepsilon\in\mathcal{A}_f}u(\varepsilon)\leq2n
\end{equation}
because each arc of $T^+$ incident to $o$ and $h_i$ where $1\leq{i}\leq{n}$ can be crossed by at most two of the arcs introduced when flipping an arc of $\mathcal{A}_f$ within $T^-$. Moreover, as $\mathcal{A}_p$ is empty, it follows from (\ref{PW.sec.4.prop.2.eq.1}) that
\begin{equation}\label{PW.sec.4.prop.2.eq.3}
\frac{n}{m}\leq|\mathcal{A}_f|\mbox{.}
\end{equation}

Combining (\ref{PW.sec.4.prop.2.eq.2}) and (\ref{PW.sec.4.prop.2.eq.3}) shows that at least one arc $\varepsilon$ in $\mathcal{A}_f$ satisfies
\begin{equation}\label{PW.sec.4.prop.2.eq.4}
u(\varepsilon)\leq2m\mbox{.}
\end{equation}

As any arc in $\mathcal{A}_f$ crosses at least $n+3m+5$ arcs of $T^+$, the number of arc crossings with $T^+$ decreases by at least $n+3m+5-u(\varepsilon)$ when $\varepsilon$ is flipped in $T^-$. By (\ref{PW.sec.4.prop.2.eq.4}) this is greater than $n+m$ as desired.
\end{proof}

According to propositions \ref{PW.sec.4.prop.1} and \ref{PW.sec.4.prop.2}, if
$$
|\mathcal{A}_f\cup\mathcal{A}_p|\geq{m+3}
$$
and
$$
2m\leq{n}\leq{m(m+3)}\mbox{,}
$$
then the arc crossings based method that estimates distances in $\mathcal{F}(\mathrm{P},X)$ will always flip an arc from $\mathcal{A}_f\cup\mathcal{A}_p$ first. We will show using Theorem~\ref{PW.sec.0.thm.0} that, for well-chosen values of $n$, $m$, and of the number of arcs in $\mathcal{A}_f\cup\mathcal{A}_p$, such a flip cannot the the first one that is performed along a geodesic from $T^-$ to $T^+$ in $\mathcal{F}(\mathrm{P},X)$. From now on, we consider a geodesic path $(T_i)_{0\leq{i}\leq{k}}$ from $T^-$ to $T^+$ in $\mathcal{F}(\mathrm{P},X)$, the blow-up triangulation $K$ of $\mathrm{P}$ associated to that path, and the corresponding family $(\Pi_i)_{0\leq{i}\leq{k}}$ of topological disks provided by Proposition \ref{PW.sec.1.lem.1}. In order to apply Theorem~\ref{PW.sec.0.thm.0} to $K$, we first prove the existence of certain arcs in this blow-up triangulation.

\begin{lem}\label{PW.sec.4.lem.1}
If $\mathcal{A}_f$ contains more than $6m+4$ arcs, then some arc in $K$ is incident to $o$ and $f$, and another one to $o$ and $p$.
\end{lem}
\begin{proof}
Let us recall that, given a vertex $x$ of $\mathrm{P}$, the distance in $\mathcal{F}(\mathrm{P},X)$ between two triangulations is at most $2|X|-2$ other than $x$ minus the total number of arcs incident to $x$ in the two triangulations \cite{SleatorTarjanThurston1988}. Taking $o$ for $x$ shows that the distance between $T^-$ and $T^+$ in $\mathcal{F}(\mathrm{P},X)$ satisfies
\begin{equation}\label{PW.sec.4.lem.1.eq.1}
k\leq2n+6m+8\mbox{.}
\end{equation}

The remainder of the proof proceeds like the proof of Lemma \ref{PW.sec.3.lem.1}. Assume that $\mathcal{A}_f$ contains more than $6m+4$ arcs. For any vertex $x$ of $\mathrm{P}$, denote by $\mathcal{S}^-_x$ the set of the tetrahedra in $K$ incident to and above an arc of $K^-$ that admits $x$ as a vertex and by $\mathcal{S}^+_x$ the set of the tetrahedra in $K$ incident to and below an arc of $K^+$ that admits $x$ as a vertex. There must be $n+1$ tetrahedra in $\mathcal{S}^-_o$ and in $\mathcal{S}^+_o$ because $o$ is incident to exactly $n+1$ arcs in $K^-\mathord{\setminus}K^+$ and to exactly $n+1$ arcs in $K^+\mathord{\setminus}K^-$. Similarly, $\mathcal{S}^+_g$ contains $n+2$ tetrahedra. Observe that $\mathcal{S}^+_o$ is disjoint from $\mathcal{S}^+_g$ because all the tetrahedra in these sets are incident to and below exactly one arc from $K^+\mathord{\setminus}K^-$. It is also disjoint from $\mathcal{S}^-_o$ because two arcs of $K$ incident to the same vertex cannot be below or above one another. Hence, if $\mathcal{S}^-_o$ and $\mathcal{S}^+_g$ were also disjoint, then $k$ would be at least $3n+4$. However, as $\mathcal{A}_f$ contains more than $6m+4$ arcs, $n$ is greater than $6m+4$, and $3n+4$ greater than $2n+6m+8$. Therefore, $k$ cannot be at least $3n+4$ because this would contradict (\ref{PW.sec.4.lem.1.eq.1}). As a consequence, some tetrahedron belongs to $\mathcal{S}^-_o\cap\mathcal{S}^+_g$. That tetrahedron has an edge with vertices $g$ and $o$ as desired. Finally, as a tetrahedron of $K$ cannot be incident to and above two arcs of $K^-$,
\begin{equation}\label{PW.sec.4.lem.1.eq.2}
|\mathcal{S}_f^-|=|\mathcal{A}_f|+2\mbox{.}
\end{equation}

By a similar argument as above, $\mathcal{S}_o^-$ is disjoint from $\mathcal{S}_f^-$ and as already mentioned, $\mathcal{S}_o^-$ is also disjoint from $\mathcal{S}_o^+$. Recall that $\mathcal{S}_o^-$ and $\mathcal{S}_o^+$ each contain $n+1$ tetrahedra. Hence, if $\mathcal{S}_o^+$ and $\mathcal{S}_f^-$ were disjoint, (\ref{PW.sec.4.lem.1.eq.2}) would implies
$$
k\geq2n+|\mathcal{A}_f|+4\mbox{.}
$$

However, as $\mathcal{A}_f$ contains more than $6m+4$ arcs, this would contradict~(\ref{PW.sec.4.lem.1.eq.1}) again. As a consequence, $\mathcal{S}^+_o\cap\mathcal{S}^-_f$ contains a tetrahedron, and that tetrahedron has an edge with vertices $f$ and $o$ as desired.
\end{proof}

\begin{figure}
\begin{centering}
\includegraphics[scale=1]{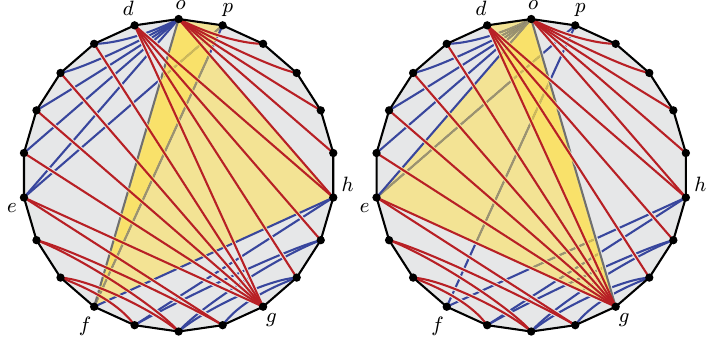}
\caption{Four triangles, two of whose are incident to the arc $\beta^-$ (left) and the others two to the arc $\beta^+$ (right).
%The four triangles $\Delta_d$, $\Delta_e$, $\Delta_h$ and $\Delta_p$ colored yellow and the arcs $\beta^-$ and $\beta^+$, colored gray. Two of these triangles are incident to $\beta^-$ (left) and two to $\beta^+$ (right).
}\label{PW.sec.4.fig.3}
\end{centering}
\end{figure}

\begin{figure}
\begin{centering}
\includegraphics[scale=1]{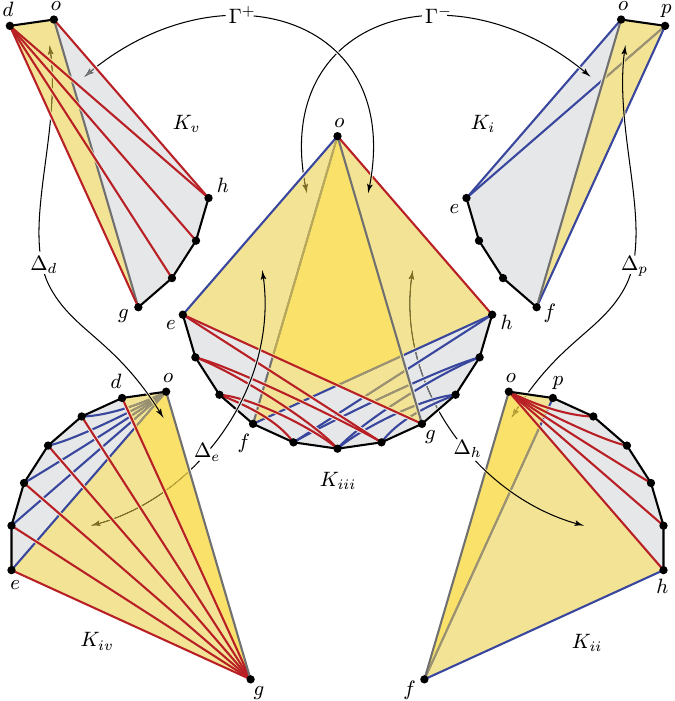}
\caption{A decomposition of $K$.}\label{PW.sec.4.fig.4}
\end{centering}
\end{figure}

Let us assume for a moment that $|\mathcal{A}_f|$ is greater than $6m+4$. By Lemma~\ref{PW.sec.4.lem.1}, some arc $\beta^-$ in $K$ has vertices $f$ and $o$ and an arc $\beta^+$ has vertices $g$ and $o$. These arcs are shown in Fig. \ref{PW.sec.4.fig.3} above $T^-$ and below $T^+$. Observe that the four triangles incident to $\beta^-$ or $\beta^+$ in the figure are each bounded by three arcs of $K$, and by Theorem \ref{PW.sec.0.thm.0}, they must belong to $K$. We will denote these four triangles by $\Delta_d$, $\Delta_e$, $\Delta_h$, and $\Delta_p$ in such a way that $\Delta_x$ admits $x$ as a vertex. As we did in Section \ref{PW.sec.3}, we will now decompose $K$ into five smaller triangulations $K_i$, $K_{ii}$, $K_{iii}$, $K_{iv}$, and $K_v$. Consider two integers $i$ and $j$ such that $\Delta_h$ is contained in $\Pi_i$ and $\Delta_e$ in $\Pi_j$. As can be seen in Fig. \ref{PW.sec.4.fig.3}, $\Delta_h$ is below $\Delta_e$ and therefore, $i$ is less than $j$. In this case, as can be seen in Fig. \ref{PW.sec.4.fig.3}, $\Pi_i$ must contain the arc $\delta$ of $K^-$ with vertices $e$ and $o$. Denote by $\Gamma^-$ the portion of $\Pi_i$ made of the points $x$ such that $\pi(x)$ is not on the same side of $\pi(\beta^-)$ as $h$ and not on the same side of $\pi(\delta)$ as $d$. By construction, $\Gamma^-$ is the union of a subset of the faces of $K$. Similarly, $\Pi_j$ must contain the arc $\varepsilon$ of $K^+$ with vertices $h$ and $o$. Let $\Gamma^+$ be the portion of $\Pi_j$ made of the points $x$ such that $\pi(x)$ is not on the same side of $\pi(\beta^+)$ as $e$ and not on the same side of $\pi(\varepsilon)$ as $p$. Again, $\Gamma^+$ is the union of a subset of the faces of $K$. Hence, cutting $K$ along $\Gamma^-$, $\Gamma^+$, and the four triangles $\Delta_d$, $\Delta_e$, $\Delta_h$, and $\Delta_p$ results in the five blow-up triangulations $K_i$, $K_{ii}$, $K_{iii}$, $K_{iv}$, and $K_v$ shown in Fig. \ref{PW.sec.4.fig.4}. The arcs contained in $\Gamma^-$ and $\Gamma^+$ are not shown in the figure as we do not know what arcs of $K$ they are. The arcs of $K^-\mathord{\setminus}K^+$ incident to $e_i$ where $1\leq{i}\leq{m}$ or to $h_i$ where $1\leq{i}\leq{n}$ are not shown either because they depend on our choice for $T^-$.

This decomposition of $K$ will allow us to determine the exact value of $k$. In order to do that it suffices to count the number of tetrahedra in each of the blow-up triangulations of that decomposition. We begin with $K_{iii}$. In the following, $l^-$ denotes the number of arcs in $K_{iii}^-\mathord{\setminus}K_{iii}^+$ incident to $o$.

\begin{lem}\label{PW.sec.4.lem.2}
There are at least $5m+5-l^-$ tetrahedra in $K_{iii}$.
\end{lem}
\begin{proof}
Consider the convex hull of the vertices of $K_{iii}$. Cutting that polygon along $\pi(\beta^+)$ results into two smaller polygons. Let us denote these smaller polygons by $\Omega$ and $\Xi$ with the convention that $\Omega$ contains $e$, and $\Xi$ contains $h$. One can see on the sketch of $K_{iii}$ at the center of Fig.~\ref{PW.sec.4.fig.4} that $K_{iii}$ satisfies the conditions of Theorem~\ref{PW.sec.1.1.thm.1} with respect to $\Xi$ and $\Omega$ when $\pi(\beta^+)$ is chosen for $\delta$. If in addition, we pick $g$ for the point $x$ in the statement of Theorem~\ref{PW.sec.1.1.thm.1}, then the blow-up triangulation $L$ provided by that theorem is the one we sketch on the left of Fig. \ref{PW.sec.4.fig.5}. It is obtained by removing from $K_{iii}$ the faces that are incident to at least two vertices of $\Xi$ and by pulling to $g$ the remaining faces that are incident to a vertex of $\Xi$ other than $g$. According to Theorem~\ref{PW.sec.1.1.thm.1}, the number of tetrahedra in $L$ is less than the number of tetrahedra in $K_{iii}$ by at least $2m+3$.

The remainder of the proof proceeds just as that of Lemma \ref{PW.sec.3.lem.2}. One should note here that the number of arcs in $L^-\mathord{\setminus}L^+$ incident to $g$ is equal to $l^-+m$. Consider the triangle in $L^-$ incident to $e$ and $g$ and denote by $\delta$ the edge of this triangle that does not admit $g$ as a vertex. Let $x$ be the vertex of $\delta$ distinct from $e$. If $x$ is equal to $e_1$, then $\pi(\delta)$ is an edge of the convex polygon $\Lambda$ whose vertices are the vertices of $L$. In that case, we denote by $\Omega$ that polygon and we denote by $\Xi$ the arc $\pi(\delta)$. 
\begin{figure}
\begin{centering}
\includegraphics[scale=1]{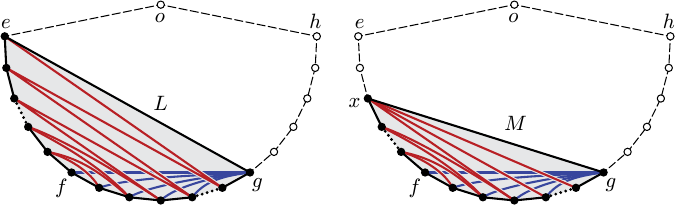}
\caption{Two blow-up triangulations considered in the proof of Lemma \ref{PW.sec.4.lem.2}. The blow-up triangulation $K_{iii}$ which they are obtained from, is sketched with dashed lines.}\label{PW.sec.4.fig.5}
\end{centering}
\end{figure}
Otherwise, we denote by $\Omega$ and $\Xi$ the two convex polygons obtained by cutting $\Lambda$ along $\pi(\delta)$ with the convention that $g$ belongs to $\Omega$ and $e_1$ to $\Xi$. One can see (for instance by looking at Fig. \ref{PW.sec.4.fig.5}) that $L$ satisfies the conditions of Theorem~\ref{PW.sec.1.1.thm.1} with respect to $\delta$, $\Xi$, and $\Omega$. Note in particular that, while $\delta$ belongs to $L^-$, the theorem can still be used, by symmetry. By Theorem~\ref{PW.sec.1.1.thm.1}, one can build, from $L$ a blow-up triangulation $M$ whose number of tetrahedra in is less than the number of tetrahedra in $L$ by at least $2|X\cap\Xi|-3$. That blow-up triangulation, sketched on the right of Fig.~\ref{PW.sec.4.fig.5}, is obtained by removing from $L$ the faces that are incident to at least two points in $X\cap\Xi$ and by pulling to $x$ the remaining faces that are incident to a point in $(X\cap\Xi)$ other than $x$. Now observe that if $x$ is not equal to $f$, then we can invoke Theorem~\ref{PW.sec.1.1.thm.1} again, by considering the triangle of $M^-$ incident to $g$ and $x$ and its edge $\delta$ that does not admit $g$ as a vertex. Repeating this procedure will eventually result in the blow-up triangulation $N$ that we already considered in the proof of Lemma \ref{PW.sec.3.lem.2} and shown at the bottom of Fig. \ref{PW.sec.3.fig.6}, on the right.

By Theorem~\ref{PW.sec.1.1.thm.1}, the number of tetrahedra in $N$ is less than the number of tetrahedra in $L$ by at least twice the number of vertices of $L$ that are no longer vertices of $N$ minus the number of times we invoked the theorem in order to change $L$ into $N$. Now recall that $L^-\mathord{\setminus}L^+$ contains $l^-+m$ arcs incident to $g$. As an arc of $L^-\mathord{\setminus}L^+$ incident to $g$ is removed each time we invoke Theorem~\ref{PW.sec.1.1.thm.1} and exactly $m$ such arcs remains in $N$, the number of tetrahedra in $N$ is less than the number of tetrahedra in $M$ by at least $2m+2-l^-$. Finally, recall that $N$ contains at least $m$ tetrahedra (because there are $m$ arcs in $N^-\mathord{\setminus}N^+$). Summing $2m+3$, $2m+2-l^-$, and $m$ results in the desired bound on the number of tetrahedra of $K_{iii}$. 
\end{proof}

When $\mathcal{A}_f$ is large enough, we can now determine the exact value of $k$.% as a consequence of Lemma \ref{PW.sec.4.lem.2}. 

\begin{lem}\label{PW.sec.4.lem.3}
If $|\mathcal{A}_f|$ is greater than $6m+4$, then $k=2n+6m+8$.
\end{lem}
\begin{proof}
Assume that $\mathcal{A}_f$ contains more than $6m+4$ arcs. We have already shown while proving Lemma \ref{PW.sec.4.lem.1} that $k$ is at most $2n+6m+8$. Hence, it suffices to establish the matching lower bound. As $\mathcal{A}_f$ contains more than $6m+4$ arcs, we know that $K$ must contain the arcs $\beta^-$, $\beta^+$ and the four triangles $\Delta_d$, $\Delta_e$, $\Delta_h$, and $\Delta_p$ shown in Fig. \ref{PW.sec.4.fig.3}. Therefore, as discussed above, $K$ can be decomposed as shown in Fig. \ref{PW.sec.4.fig.4} into five blow-up triangulations $K_i$, $K_{ii}$, $K_{iii}$, $K_{iv}$, and $K_v$. Recall that this decomposition requires one to cut $K$ along $\Gamma^-$ and $\Gamma^+$ and these four triangles.

Now let us count the number of tetrahedra in each piece of the decomposition. According to Lemma \ref{PW.sec.4.lem.2}, $K_{iii}$ contains at least $5m+5-l^-$ tetrahedra, where $l^-$ is the numbers of arcs incident to $o$ contained in $K_{iii}^-\mathord{\setminus}K_{iii}^+$. Let us count the number of tetrahedra in the other four blow-up triangulations. Observe that $K_v$ contains at least $m+1$ tetrahedra because $K_v^-\mathord{\setminus}K_v^+$ contains exactly $m+1$ arcs. The situation is different for $K_i$ because there may be less than $m+1$ arcs in $K_i^-\mathord{\setminus}K_i^+$. However, we know that $o$ is incident to $l^-$ arcs from $K_i^+\mathord{\setminus}K_i^-$ because $K_i$ and $K_{iii}$ are glued along $\Gamma^-$ (see Fig. \ref{PW.sec.4.fig.4}). In particular, none of these arcs can be contained in both $K_i^-$ and $K_i^+$ because they are above the arc in $K_i^-$ with vertices $e$ and $p$. As a consequence, the number of tetrahedra in $K_i$ is at least $l^-$. Finally, observe that $K_{ii}$ and $K_{iv}$ must each contain at least $n+1$ tetrahedra. Summing $5m+5-l^-$, $m+1$, $l^-$, and two times $n+1$ provides the desired bound.
\end{proof}

Interestingly, Lemma \ref{PW.sec.4.lem.3} shows that the number of tetrahedra in $K$, or equivalently the distance of $T^-$ and $T^+$ in $\mathcal{F}(\mathrm{P},X)$ does not depend on $\mathcal{A}_f$ or $\mathcal{A}_p$ except maybe when $\mathcal{A}_f$ contains few arcs. In particular, when
$$
|\mathcal{A}_f|>6m+4
$$
and
$$
2m\leq{n}\leq{m(7m+5)}\mbox{,}
$$
it follows from Lemma \ref{PW.sec.4.lem.3}, Proposition \ref{PW.sec.4.prop.1}, and Proposition \ref{PW.sec.4.prop.2} that the first flip performed by the arc crossings-based distance estimation method will not be the first flip in any geodesic path from $T^-$ to $T^+$ in $\mathcal{F}(\mathrm{P},X)$. In fact, we have proven that the triangulation resulting from that flip will be at the same distance of $T^+$ in $\mathcal{F}(\mathrm{P},X)$ than $T^-$.

\begin{thm}\label{PW.sec.4.thm.1}
Assume that $m$ is positive and $n$ is equal to $m(7m+5)$. If $T^-$ contains the comb triangulations of $\mathrm{P}_h$ at $f$ and of $\mathrm{P}_e$ at $p$, then
$$
\frac{\ell(T^-,T^+)}{d(T^-,T^+)}\geq\frac{3}{2}-\frac{9m+9}{14m^2+16m+8}\mbox{.}
$$
\end{thm}
\begin{proof}
Let us assume that $T^-$ contains, as a subset, the comb triangulation of $\mathrm{P}_h$ at $f$ and the comb triangulation of $\mathrm{P}_e$ at $p$. Observe that by our choices for $n$, $m$, and $T^-$, the assumptions of Proposition~\ref{PW.sec.4.prop.1}, Proposition~\ref{PW.sec.4.prop.2}, and Lemma~\ref{PW.sec.4.lem.3} are satisfied for $T^-$ and all the triangulations obtained from it by flipping at most $n-6m-5$ of the arcs contained in $\mathcal{A}_f\cup\mathcal{A}_p$. In particular, the distance in $\mathcal{F}(\mathrm{P},X)$ of any of these triangulations to $T^+$ is always $2n+6m+8$. Moreover, by Propositions~\ref{PW.sec.4.prop.1} and~\ref{PW.sec.4.prop.2}, the first $n-6m-5$ flips performed by the arc crossings-based method will only flip arcs of $T^-$ contained in $\mathcal{A}_f\cup\mathcal{A}_p$. As a consequence, the method overestimates the distance of $T^-$ and $T^+$ by a factor of at least
$$
1+\frac{n-6m-5}{2n+6m+8}\mbox{.}
$$

Replacing $n$ by $m(7m+5)$ in this expression completes the proof.
\end{proof}

\begin{figure}
\begin{centering}
\includegraphics[scale=1]{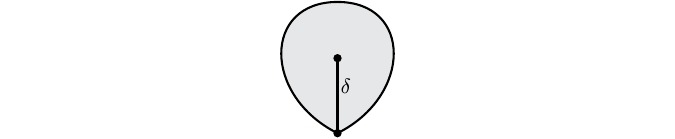}
\caption{A self-folded triangle.}\label{PW.sec.4.fig.6}
\end{centering}
\end{figure}

Let us finally consider the case of an arbitrary orientable topological surface $\Sigma$. We assume that $\Sigma$ has finite genus and a finite number of boundary components. The set $X$ that serves as the vertex set for the triangulations is made of punctures of $\Sigma$ and marked points on the boundary of $\Sigma$ in such a way that at least one marked point is contained in each boundary component of $\Sigma$. In this context, a triangulation of $\Sigma$ with vertex set $X$ is a set $T$ of pairwise non-crossing and pairwise non-isotopic arcs between points in $X$ that is maximal with respect to inclusion. The arcs of $T$ are considered up to proper isotopy. In other words, they are isotopy classes of non-oriented paths that connect two vertices of $T$, and we ask that these classes can be realized disjointly. Some arcs of $T$ may be loops (when their ends coincide), and several arcs may connect the same pair of vertices, but $T$ is still a triangulation in the sense that the complement in $\Sigma$ of the union of its arcs is a finite collection of triangles. As in the case of a polygon, a flip removes an arc from $T$ that is incident to two distinct triangles and replaces it with the only other arc such that the resulting set is a triangulation of $\Sigma$. It should be noted that not all arcs can be flipped. For instance, the arc $\delta$ shown in Fig.~\ref{PW.sec.4.fig.6} is not flippable because it is not incident to two different triangles. %Apart from arcs that belong, up to isotopy, to the boundary of the surface, this is the only example of a non-flippable arc.
The flip-graph $\mathcal{F}(\Sigma,X)$ of $\Sigma$ can then be defined just the same way as for polygons. This setup is considered for instance in \cite{DisarloParlier2019}.

Now observe that Theorem~\ref{PW.sec.0.thm.3} immediately follows from Theorem~\ref{PW.sec.4.thm.1} in the special case when $\Sigma$ is a disk. Indeed, the flip-graph of a disk with $n$ marked points in its boundary and no puncture is isomorphic to the flip-graph of a convex Euclidean polygon with $n$ vertices. Let us prove that theorem in general as a consequence of Theorem~\ref{PW.sec.4.thm.1}.

\begin{proof}[Proof of Theorem~\ref{PW.sec.0.thm.3}]
Let us place in $X$ marked points from the boundary of $\Sigma$ and punctures of $\Sigma$, making sure that there is at least one marked point in each boundary component. Consider a closed curve $\gamma$ through all the points from $X$ in such a way that one of the surfaces obtained when $\Sigma$ is cut along $\gamma$ is a disk $\Delta$ bounded by $\gamma$. In particular, $\mathcal{F}(\Delta,X)$ is isomorphic to the flip-graph of a convex polygon with $|X|$ vertices. Now, make $X$ large enough so that we can embed the triangulations $T^-$ and $T^+$ from Theorem~\ref{PW.sec.4.thm.1} as triangulations of $\Delta$ with vertex set $X$ and complete these triangulations into two triangulations $U^-$ and $U^+$ of $\Sigma$ by adding the same set $T$ of arcs. As, for any arc $\varepsilon$ between two points of $X$, the subgraph $\mathcal{F}_\varepsilon(\Sigma,X)$ is strongly convex in $\mathcal{F}(\Sigma,X)$ \cite{DisarloParlier2019}, the distance of $U^-$ and $U^+$ in $\mathcal{F}(\Sigma,X)$ will be the same as in the distance of $T^-$ and $T^+$. Moreover, the crossings-based distance estimation method will perform the same sequence of flips as from $T^-$ to $T^+$ and the result follows from Theorem \ref{PW.sec.4.thm.1}
\end{proof}

\medskip
\noindent{\bf Acknowledgement.} Both authors are partially supported by the ANR project SoS (Structures on Surfaces), grant number ANR-17-CE40-0033. The second author is supported by the National Natural Science Foundation of China (Grant No. 12301432) and the Natural Science Foundation of Guangdong Province (Grant No. 2023A1515010658). The manuscript was reworked while the second author was visiting the first one thanks to funding of the MathSTIC (CNRS FR3734) research consortium. We are grateful to the referees for carefully reading the manuscript and providing comments and suggestions that greatly improved the exposition.

%\medskip
%\noindent{\bf Statements.} The authors state that there is no conflict of interest. This article has no associated data.

\bibliography{StrongConvexity}
\bibliographystyle{ijmart}

\end{document}